\documentclass[pdflatex,sn-mathphys-num]{sn-jnl}

\usepackage{graphicx}%
\usepackage{multirow}%
\usepackage{amsmath,amssymb,amsfonts}%
\usepackage{amsthm}%
\usepackage{mathrsfs}%
\usepackage[title]{appendix}%
\usepackage{xcolor}%
\usepackage{textcomp}%
\usepackage{manyfoot}%
\usepackage{booktabs}%
\usepackage{algorithm}%
\usepackage{algorithmicx}%
\usepackage{algpseudocode}%
\usepackage{listings}%
\usepackage{mathtools}
\usepackage{enumitem}

\theoremstyle{thmstyleone}%
\newtheorem{theorem}{Theorem}[section]%
\newtheorem{proposition}[theorem]{Proposition}%
\newtheorem{corollary}[theorem]{Corollary}%
\newtheorem{lemma}[theorem]{Lemma}%
\newtheorem{remark}[theorem]{Remark}%
\newtheorem{definition}[theorem]{Definition}%

\newcommand{\ZZ}{\mathbb{Z}}
\newcommand{\RR}{\mathbb{R}}

\newcommand{\PP}{\mathbb{P}}
\newcommand{\EE}{\mathbb{E}}
\newcommand{\TT}{\mathbb{T}}
\newcommand{\BB}{\mathbb{B}}
\newcommand{\var}{\mathrm{Var}}

\newcommand{\sC}{\mathcal{C}}

\newcommand{\sH}{\mathcal{H}}

\newcommand{\sN}{\mathcal{N}}

\newcommand{\diff}{{\mathrm d}}
\newcommand{\dist}{\mathrm{dist}}
\newcommand{\ds}{\diff s}

\newcommand{\dx}{\diff x}
\newcommand{\dt}{\diff t}

\newcommand{\ind}{\mathbf{1}}

\newcommand{\ig}{\mathrm{L}}
\newcommand{\igxi}{\mathrm{L}^{\xi}}
\newcommand{\igl}{\mathrm{L}_{\boxn, \beta}^{\xi}}

\newcommand{\gibbsl}{\rho_{\boxn, \beta}^{\xi}}

\newcommand{\gamxi}{\Gamma^{\xi}}
\newcommand{\gaml}{\Gamma_{\boxn, \beta}^{\xi}}

\newcommand{\torusn}{\TT_{n}}
\newcommand{\torusnd}{\TT_{n}^{d}}
\newcommand{\boxn}{\Lambda_{n}}

\newcommand{\s}{\sigma}

\newcommand{\sy}{\sigma^{y}}

\newcommand{\delx}{\nabla_{x}}
\newcommand{\dely}{\nabla_{y}}

\newcommand{\avg}[1]{\left\langle #1 \right\rangle}

\newcommand{\avgl}[1]{\left\langle #1 \right\rangle_{\boxn, \beta}^{\xi}}

\newcommand{\interior}[1]{{\kern0pt#1}^{\mathrm{o}}}

\newcommand{\<}{\langle}
\renewcommand{\>}{\rangle}

\newcommand{\distk}{d_{\mathrm{K}}}

\begin{document}

\title[Article Title]{Berry--Esseen Bounds for High-Temperature Ising Models}


\author[1]{\fnm{Timothy M.} \sur{Garoni}}\email{tim.garoni@monash.edu}
\author[1]{\fnm{Aram} \sur{Perez}}\email{aram.perez@monash.edu}

\affil[1]{\orgdiv{School of Mathematics}, \orgname{Monash University}, \orgaddress{\city{Clayton}, \postcode{3800}, \state{VIC},
    \country{Australia}}}

\abstract{We consider the nearest-neighbour Ising model on \(\ZZ^d\), and establish sharp results for the normal approximation of the magnetisation throughout the high-temperature regime. In particular, we establish optimal Berry–Esseen bounds between the magnetisation in an arbitrary finite region and a Gaussian with the corresponding finite-volume variance. These results hold for finite-volume Gibbs measures with arbitrary boundary conditions, as well as for the infinite-volume Gibbs measure. We also determine sharp rates of convergence of the magnetisation in a finite box to a Gaussian with variance given by the infinite-volume susceptibility. For this fixed Gaussian target, we show that, in dimensions \(d\ge 3\), periodic boundary conditions yield strictly faster convergence than free boundary conditions. Our approach uses Stein’s method to quantify the normal approximation error by analysing the semigroup of the continuous-time Ising heat-bath dynamics.}




\maketitle

\section{Introduction}
Limit theorems for the total magnetisation of spin systems have been well studied for over half a century. In particular, Newman~\cite{Newman1980} established general central limit theorems for translation-invariant FKG systems with finite susceptibility. These theorems apply, in particular, to Ising models throughout
the high-temperature regime. A natural question associated to these limit theorems is: what are their corresponding rates of convergence? It may be natural to expect that, for a system of weakly dependent random variables, one may recover the Berry--Esseen rate obtained for the classical central limit theorem for sums of i.i.d. random variables.

Stein's method \cite{Stein1972, Stein1986} provides a set of techniques for establishing explicit rates of convergence for limit theorems, which is particularly powerful when analysing sums of dependent random variables. For normal approximation, it is centred around the observation that a random variable $Z$ has law $\mathcal{N}(0, \sigma^{2})$ if and only if it satisfies the integration by parts identity
\begin{equation}\label{eq: stein identity}
    \sigma^{2}\EE[f'(Z)] = \EE[Zf(Z)], \quad  f\in \sC,
\end{equation}
where $\sC$ is a sufficiently large class of test functions. For a mean-zero random variable $X$, Stein's method bounds distances between the laws of $X$ and $Z\sim \mathcal{N}(0, \sigma^{2})$ by quantifying how close $X$ is to satisfying this identity.

Stein's theory has been successfully used to establish convergence rates for a range of mean-field spin models in both non-critical and critical regimes, including the Curie--Weiss model~\cite{chatterjee2011nonnormal, EichelsbacherLowe2010}, the Curie--Weiss--Potts model \cite{EichelsbacherMartschink2015}, the $O(N)$ model~\cite{KirkpatrickNawaz2016, Garoni2026} and the Blume--Emery--Griffiths model~\cite{eichelsbacher2014rates}.
These works employ the exchangeable-pairs approach to Stein’s method, analysing conditional expectations associated with a suitable coupling of the total magnetisation. 
For spin systems, a natural construction is provided by taking a single step in the Glauber (heat-bath) dynamics at equilibrium, where one updates the spin at a uniformly random vertex according to the conditional distribution obtained by fixing the remaining spins. The lack of spatial geometry in the mean-field setting makes the computation of these conditional expectations rather simple, since the conditional expectations essentially depend on the total magnetisation itself. This makes bounding the errors required to establish explicit convergence rates sufficiently tractable.

For lattice models, results on convergence rates are relatively sparse and typically lead to suboptimal rates of convergence. The method of exchangeable pairs using the Glauber construction has yet to be applied successfully due to the local dependence on the updated spin, making the conditional expectations more difficult to control. In the case of the Ising model, other methods, including Fourier estimates~\cite{Bulinskii1996}, martingale embeddings~\cite{fang2025}, and a Stein's method approach adapted to FKG systems~\cite{Goldstein2018}, have yet to achieve optimal convergence rates throughout the entire high-temperature regime. 

In this work, we consider high-temperature Ising models on the integer lattice, where the total magnetisation has asymptotically normal fluctuations. Rather than relying on the exchangeable pair approach, we develop a method based on ideas from the theory of Stein kernels and Stein couplings; see, e.g., \cite{Ley2017, Courtrade2019} for an overview of Stein kernels and \cite{chen2010steincouplingsnormalapproximation} for Stein couplings. This reduces  non-asymptotic estimates for the distance between the laws of the total magnetisation and a corresponding normal distribution to variance estimates for an associated random variable constructed from the Glauber semigroup. Using this approach, we establish optimal Berry--Esseen bounds throughout the entire subcritical regime. We also compare the convergence rates for finite-volume Ising models under free and periodic boundary conditions, showing that, when $d\geq 3$, periodic boundary conditions yield strictly faster convergence to the normal distribution with asymptotic variance given by the infinite-volume susceptibility.

\subsection{The Ising model}  
Throughout, $d\ge1$ will denote a fixed integer. For each integer $n\geq 3$, let $\boxn:= [-n/2, n/2)^{d}\cap \ZZ^{d}$. Define the graphs $\BB_{n}^{d} := (\boxn, E_{\BB}(\boxn))$ and $\torusnd := (\boxn, E_{\TT}(\boxn))$ where
\begin{equation}\label{eq: EB def}
    E_{\BB}(\Lambda) :=\{\{x, y\}: x-y=\pm e_{i},\, x, y\in \Lambda,\, i\in [d]\}, \quad \Lambda \subseteq \ZZ^{d},
\end{equation}
and
\begin{equation}
    E_{\TT}(\Lambda_{n}) :=\{\{x, y\}:x-y=\pm e_{i} \bmod n,\, x, y\in \Lambda_{n},\, i\in [d]\}.
  \end{equation}
  Here $[d]\coloneq\{1,2,\ldots,d\}$, and $e_i$ denotes the $i$th standard basis vector in $\RR^d$.
Let $E(\ZZ^{d}):=E_{\BB}(\ZZ^{d})$ denote the nearest-neighbour edges induced by the vertices in $\ZZ^{d}$ and, as short-hand, abbreviate an edge $\{x, y\}$ by $xy$. We consider the Ising model on $\BB_{n}^{d}$, $\torusnd$ and $\ZZ^{d}$.

Let $\Omega := \{-1, 1\}$ denote the single-spin space for the Ising model. For any nonempty finite $V\subset\ZZ^d$,
the Gibbs measure for the Ising model on $(V,E_{\BB}(V))$ is defined at inverse temperature $\beta\geq0$ and with boundary condition $\xi \in \{-1, 0, 1\}^{\ZZ^{d}}$ via
\begin{equation}\label{Eq: gibbs measure}
    \rho_{V, \beta}^{\xi}(\sigma) := \frac{1}{Z_{V, \beta}^{\xi}}\exp(-\beta \sH_V^{\xi}(\sigma)), \quad \sigma \in \Omega^{V}
\end{equation}
where
\begin{equation}\label{Eq: hamiltonian}
  \sH_V^{\xi}(\sigma) := -\sum_{\substack{xy\in E(\ZZ^{d})\\ x, y \in V}}\sigma(x)\sigma(y)
  -\sum_{\substack{xy\in E(\ZZ^{d})\\x\in V, y\notin V}}\sigma(x)\xi(y).
\end{equation}
is the Hamiltonian and
\begin{equation}\label{Eq: partition function}
    Z_{V, \beta}^{\xi} := \sum_{\sigma \in \Omega^{V}}\exp(-\beta \sH_V^{\xi}(\sigma))
\end{equation}
is the partition function. We will almost exclusively be concerned with the case $V=\boxn$.
Note that the case $\xi(y) = 0, \ \ \forall  y\in \ZZ^{d}$, denoted $\xi = 0$ henceforth, corresponds to free boundary conditions, whilst the case $\xi(y) = +1, \ \ \forall  y\in \ZZ^{d}$, denoted by $\xi = +$, corresponds to plus boundary conditions.

In addition, we consider the model on the torus $\torusnd = (\boxn, E_{\TT}(\boxn))$, defined by replacing the Hamiltonian in~\eqref{Eq: hamiltonian} by
\begin{equation}\label{eq: hamiltonian torus}
    \sH^{\mathrm{per}}(\sigma) := -\sum_{xy\in E_{\TT}(\boxn)}\sigma(x)\sigma(y).
\end{equation}
We let $\rho_{\boxn, \beta}^{\mathrm{per}}$ and $Z_{\boxn, \beta}^{\mathrm{per}}$ denote the corresponding Gibbs measure and partition function for the torus model. Henceforth, we write $\xi = \mathrm{per}$ when considering the Ising model on the torus. 

Let $\langle  \, \cdot \, \rangle_{\boxn, \beta}^{\xi}$ denote the expectation with respect to the Gibbs measure $\rho_{\boxn, \beta}^{\xi}$. For a given configuration $\sigma$ and set $\Lambda \subseteq \boxn$, we denote the total magnetisation in $\Lambda$ by $S_{\Lambda}=S_{\Lambda}(\sigma) := \sum_{x\in \Lambda}\sigma(x)$
and the magnetisation density in $\Lambda$ by $m_{\Lambda} =m_{\Lambda}(\sigma) := S_{\Lambda}(\sigma)/|\Lambda|$. We will often refer to these quantities as simply the total magnetisation and the magnetisation density when $\Lambda$ is clear from context. For given $\beta\ge0$ and boundary condition $\xi$, we additionally define 
\begin{equation}
    W_{\Lambda}=W_{\Lambda}(\sigma) := \frac{1}{\sqrt{|\Lambda|}}\left(S_{\Lambda} - \avgl{S_{\Lambda}}\right)\big.
\end{equation}
The finite-volume susceptibility corresponding to $\Lambda\subseteq \boxn$ under the Gibbs measure $\rho_{\boxn, \beta}^{\xi}$ is then defined as
\begin{equation}\label{eq: sus def}
    \chi_{\boxn, \beta}^{\xi}(\Lambda) := \avgl{W_{\Lambda}^{2}} = \frac{1}{|\Lambda|}\sum_{x, y \in \Lambda}\left(\avgl{\sigma(x) \sigma(y)} - \avgl{\sigma(x)}\avgl{\sigma(y)}\right).
\end{equation}

A function $f:\Omega^{\ZZ^{d}}\to \RR$ is said to be \textit{local} if there exists a finite set $S\subset \ZZ^{d}$ such that $f$ depends only on $(\sigma(x))_{x\in S}$.
Following \cite[Definition 3.14]{FriedliVelenik2017}, we then say that the sequence of finite-volume Gibbs measures $\left(\rho_{\boxn, \beta}^{\xi}\right)_{n\in \ZZ^{+}}$ converges to the infinite-volume state $\avg{\,\cdot\,}_{\ZZ^{d}, \beta}^{\xi}$, which we denote by $\lim_{n\to \infty}\langle \, \cdot \,\rangle_{\boxn, \beta}^{\xi} = \langle  \,\cdot\, \rangle_{\ZZ^{d}, \beta}^{\xi}$, if
\begin{equation}
    \lim_{n\to \infty}\langle  f \rangle_{\boxn, \beta}^{\xi} = \langle  f \rangle_{\ZZ^{d}, \beta}^{\xi}
\end{equation}
for every local function $f:\Omega^{\ZZ^{d}} \to \RR$. Equivalently, the sequence of finite-volume Gibbs measures converges weakly to a corresponding infinite-volume Gibbs measure $\rho_{\ZZ^{d}, \beta}^{\xi}$ (see, e.g., \cite{Georgii2011} for a discussion). 

Let
\begin{equation}
    m^{*}(\beta):=\<\sigma(0)\>_{\ZZ^{d}, \beta}^{+} = \lim_{n \to \infty}\<\sigma(0)\>_{\boxn, \beta}^{+},
\end{equation}
be the spontaneous magnetisation. By the GKS inequalities~\cite{Griffiths1967II, KellySherman1968} the spontaneous magnetisation is non-decreasing, which implies that
\begin{equation}\label{eq: betac}
    \beta_{c}:=\inf\{\beta>0:m^{*}(\beta)>0\} = \sup\{\beta>0:m^{*}(\beta) = 0\}.
\end{equation}
Due to \cite[Theorem 2]{Lebowitz1972}, there is a unique infinite-volume state independent of the boundary condition if and only if $m^{*}(\beta) = 0$. Consequently, for all $\beta < \beta_{c}$\footnote{When $d\geq 2$, equation~\eqref{eq: unique Gibbs state} can be extended to all $\beta\leq \beta_{c}$ thanks to the proofs that the spontaneous magnetisation vanishes at $\beta = \beta_{c}$; see, e.g., \cite{Yang1952, AizenmanDCSid2015}. When $d=1$, one simply has $\beta_c=+\infty$.}
\begin{equation}\label{eq: unique Gibbs state}
    \lim_{n\to \infty}\langle  \, \cdot\, \rangle_{\boxn, \beta}^{\xi} = \langle  \, \cdot\, \rangle_{\ZZ^{d}, \beta}, \quad \forall \,\xi\in \{-1, 0, 1\}^{\ZZ^{d}}\cup\{\mathrm{per}\}.
\end{equation}
The state $ \langle  \, \cdot\, \rangle_{\ZZ^{d}, \beta}$ is translation invariant and, moreover, satisfies $ \langle  \sigma(0)\rangle_{\ZZ^{d}, \beta} = 0$ and exponential decay of correlations in the entire high-temperature regime:
\begin{equation}\label{eq: exp decay tp}
    0 \leq \langle  \sigma(0)\sigma(x) \rangle_{\ZZ^{d}, \beta} \leq \exp(-m_{\beta, d}|x|_{1}),\quad\quad \beta<\beta_{c},
\end{equation}
where $|x|_{1}:= |x_{1}|+\cdots+|x_{d}|$ denotes the $L^{1}$ norm of $x\in \ZZ^{d}$. We may extend the definition of the finite-volume susceptibility~\eqref{eq: sus def} to the infinite-volume Gibbs measure so that for finite $\Lambda\subset \ZZ^{d}$, $W_{\Lambda} = S_{\Lambda}/\sqrt{|\Lambda|}$ and
\begin{equation}
    \chi_{\ZZ^{d}, \beta}(\Lambda) :=\avg{W_{\Lambda}^2}_{\ZZ^{d}, \beta} = \frac{1}{|\Lambda|}\sum_{x, y\in \Lambda} \avg{\sigma(x)\sigma(y)}_{\ZZ^{d}, \beta}
\end{equation}
The infinite-volume susceptibility $\chi_{\ZZ^{d}, \beta} :=\lim_{n\to \infty}\chi_{\ZZ^{d}, \beta}(\boxn)$ then satisfies
\begin{equation}\label{eq: finite sus}
    \chi_{\ZZ^{d}, \beta} = \sum_{x\in \ZZ^{d}}\langle\sigma(0) \sigma(x)\rangle_{\ZZ^{d}, \beta} < \infty,\quad\quad \beta<\beta_{c}.
\end{equation}
Exponential decay of correlations~\eqref{eq: exp decay tp} and the finiteness of the infinite-volume susceptibility~\eqref{eq: finite sus} were first established in~\cite{Aizenman1987}. 

\subsection{Statement of main results}
In this paper we calculate, for all $\beta < \beta_{c}$, Berry--Esseen bounds between the law of $W_{\Lambda}$ under $\rho_{\boxn, \beta}^{\xi}$ and a normal random variable with the same variance. Our main result is the following.
\begin{theorem}\label{Theorem: main result 1}
  Let $\beta\in[0,\beta_{c})$, and let $\Lambda \subset \ZZ^{d}$ be finite and nonempty.
  Then there are positive constants $c = c(\beta, d)$ and $C = C(\beta, d)$ such that, for every $n\ge3$ satisfying $\Lambda_{n}\supseteq \Lambda$, and every boundary condition $\xi$, 
    \begin{equation}\label{eq: berry esseen same variance}
        \frac{c}{\sqrt{|\Lambda|}} \leq \sup_{z \in \RR}\left|\rho_{\boxn, \beta}^{\xi}\left(W_{\Lambda} \leq z\right) - \PP\left(Y_{\boxn, \beta}^{\xi}(\Lambda) \leq z\right)\right| \leq \frac{C}{\sqrt{|\Lambda|}},
    \end{equation}
    where $Y_{\boxn, \beta}^{\xi}(\Lambda) \sim \sN\left(0, \chi_{\boxn, \beta}^{\xi}(\Lambda)\right)$.
\end{theorem}

Central limit theorems for the magnetisation of Ising systems on $\ZZ^{d}$ were first proved by Newman \cite{Newman1980} using an elegant blocking argument. More precisely it was shown that, under the infinite-volume Gibbs measure $\rho_{\ZZ^{d}, \beta}$, the law of $W_{\boxn}$ converges in distribution to $\mathcal{N}(0, \chi_{\ZZ^{d}, \beta})$. The method applied, more generally, to positively associated random fields on $\ZZ^{d}$, i.e. fields satisfying the FKG inequality \cite{FKG1971}. Subsequently, rates of convergence were established by Bulinskii~\cite{Bulinskii1996} for positively associated random fields, which held under mild moment assumptions and exponential decay of correlations. Applied to the high-temperature Ising model on $\ZZ^{d}$ under the infinite-volume Gibbs measure, \cite[Theorem 2]{Bulinskii1996} established Berry--Esseen bounds of order $|\Lambda|^{-1/2}\log(|\Lambda|)^{d}$ between $W_{\Lambda}$ and a normal random variable with the same variance. More recently, a general approach for establishing such bounds using positive correlations was developed in \cite{Goldstein2018} using Stein's method \cite{Stein1972, Stein1986}. Applied to the Ising model on $\ZZ^{d}$, convergence rates of order $|\Lambda|^{-1/(2d+2)}$ were obtained in the Wasserstein metric, using an analogous blocking argument to \cite{Newman1980} (see, e.g., \cite[Theorem 2.1]{Goldstein2018}). Finite-volume normal approximation results were established in~\cite{fang2025} for general Ising systems satisfying sufficiently strong high-temperature conditions, using martingale embeddings and a Poincar\'e inequality.
In particular, \cite[Theorem 1.1]{fang2025} gives bounds of order $O(|\Lambda|^{-1/14})$ in the 2-Wasserstein distance between the normalised magnetisation and a Gaussian with the same variance, for sufficiently high temperatures.

Theorem~\ref{Theorem: main result 1} is the first estimate to yield order $|\Lambda|^{-1/2}$ convergence for finite-volume Ising models, for all $\beta < \beta_{c}$. We note, however, that the generality of Theorem~\ref{Theorem: main result 1} allows us to additionally establish convergence rates under the infinite-volume Gibbs measure $\rho_{\ZZ^{d}, \beta}$. The following corollary improves the result of Bulinskii~\cite{Bulinskii1996} in the specific case of the high-temperature Ising model on $\ZZ^{d}$.
\begin{corollary}\label{corol: infinite volume limit theorem}
  Let $\beta\in[0,\beta_{c})$, let $\Lambda \subset \ZZ^{d}$ be finite and nonempty, and let $W_{\Lambda} = S_{\Lambda}/\sqrt{|\Lambda|}$.
  Then there are positive constants $c = c(\beta, d)$ and $C = C(\beta, d)$ such that
    \begin{equation}\label{eq: berry esseen same variance inf}
        \frac{c}{\sqrt{|\Lambda|}} \leq \sup_{z \in \RR}\left|\rho_{\ZZ^{d}, \beta}\left(W_{\Lambda} \leq z\right) - \PP\left(Y_{\ZZ^{d}, \beta}(\Lambda) \leq z\right)\right| \leq \frac{C}{\sqrt{|\Lambda|}},
    \end{equation}
    where $Y_{\ZZ^{d}, \beta}(\Lambda)\sim \mathcal{N}(0, \chi_{\ZZ^{d}, \beta}(\Lambda))$.
\end{corollary}

To quantify the approach to the thermodynamic limit at the level of fluctuations, a natural question is to determine the rate of convergence to a fixed Gaussian target whose variance is given by the infinite-volume susceptibility \(\chi_{\ZZ^{d},\beta}\). We consider the case \(\Lambda=\Lambda_n\). As a corollary of Theorem~\ref{Theorem: main result 1}, we establish the following.
\begin{corollary}\label{corol: stationary gaussian}
    Let $0<\beta < \beta_{c}$ and $Y_{\beta}\sim \sN\left(0, \chi_{\ZZ^{d}, \beta}\right)$. Then there are positive constants $c_{\mathrm{per}}$, $C_{\mathrm{per}}$, $c_{0}$, $C_{0}$, $ c_{\infty}$, $ C_{\infty}$ (depending only on $\beta$ and $d$) such that, for every $n\ge3$, 
    \begin{equation}\label{eq: stationary gaussian per}
        \frac{c_{\mathrm{per}}}{n^{d/2}} \leq\sup_{z \in \RR}|\rho_{\boxn, \beta}^{\mathrm{per}}(W_{\boxn} \leq z) - \PP(Y_{\beta}\leq z)| \leq \frac{C_{\mathrm{per}}}{n^{d/2}},
    \end{equation}
    \begin{equation}\label{eq: stationary gaussian free}
        \frac{c_{0}}{n \wedge n^{d/2}}\leq\sup_{z \in \RR}|\rho_{\boxn, \beta}^{0}(W_{\boxn} \leq z) - \PP(Y_{\beta}\leq z)| \leq \frac{C_{0}}{n \wedge n^{d/2}},
    \end{equation}
    and
    \begin{equation}\label{eq: stationary gaussian inf}
        \frac{c_{\infty}}{n \wedge n^{d/2}}\leq\sup_{z \in \RR}|\rho_{\ZZ^{d}, \beta}(W_{\boxn} \leq z) - \PP(Y_{\beta}\leq z)| \leq \frac{C_{\infty}}{n \wedge n^{d/2}}.
    \end{equation}
\end{corollary}  
Corollary~\ref{corol: stationary gaussian} shows that when $d\geq3$, the law of $W_{\boxn}$ under $\rho_{\boxn, \beta}^{\mathrm{per}}$ converges significantly faster to the law of $Y_{\beta}$, compared with the law of $W_{\boxn}$ under $\rho_{\boxn, \beta}^{0}$ or $\rho_{\ZZ^{d}, \beta}$. The slower convergence rate in~\eqref{eq: stationary gaussian free} (resp.~\eqref{eq: stationary gaussian inf}) is a consequence of the slower convergence of $\chi_{\boxn, \beta}^{0}(\boxn) \to \chi_{\ZZ^{d}, \beta}$ (resp. $\chi_{\ZZ^{d}, \beta}(\boxn) \to \chi_{\ZZ^{d}, \beta}$), which we show to be $\Theta(1/n)$ as $n \to \infty$; see Lemma~\ref{lemma: sus bounds free and periodic}. On the other hand $\chi_{\boxn, \beta}^{\mathrm{per}}(\boxn)\to \chi_{\ZZ^{d}, \beta}$ at a rate of $O(e^{-cn})$. Equation~\eqref{eq: stationary gaussian inf} establishes the optimal convergence rate for Newman's central limit theorem~\cite{Newman1980} in the case of the subcritical, nearest-neighbour Ising model on $\ZZ^{d}$.

Our strategy for proving Theorem~\ref{Theorem: main result 1} is to use Stein's method for normal approximation. Our implementation will rely on the study of the continuous-time heat-bath dynamics for the Ising model, which we formally introduce now.

\subsection{Heat-bath process}
We begin by providing a concrete construction for the Ising Glauber process with heat-bath updates, and then recall some of its general properties. See, e.g.,~\cite{Martinelli1999} for more details.

Let $\sigma\in\Omega^{\boxn}$ and $x\in\boxn$.
The local magnetisation at site $x$, in configuration $\sigma$, corresponding to boundary condition $\xi\in \{-1, 0, 1\}^{\ZZ^{d}}$ is 
\begin{equation}
  S_{x}^{\xi}(\sigma) := \sum_{\substack{y \in \boxn\\ xy\in E_{\BB}(\boxn)}}\sigma(y) + \sum_{\substack{y \notin \boxn\\ xy\in E(\ZZ^{d})}}\xi(y),
\end{equation}
while for periodic boundary conditions
\begin{equation}
S_{x}^{\mathrm{per}}(\sigma) := \sum_{\substack{y \in \boxn\\xy\in E_{\TT}(\boxn)}}\sigma(y).
\end{equation}
Now set, for any $\xi\in\{-1,0,1\}^{\ZZ^d}\cup\{\mathrm{per}\}$, 
\begin{equation}
  p_x^{\xi}(\sigma)\coloneqq \frac{1}{2}\left[1+\tanh\left(\beta S_{x}^{\xi}(\sigma)\right)\right].
\end{equation}
 
\begin{definition}\label{def: Glauber dynamics}
  Let $\{U_i,V_i\}_{i\in\ZZ^+}$ be an independent family of random variables, with $U_i$ uniform on $[0,1]$ and $V_i$ uniform on $\boxn$, and let $(N_t)_{t\ge0}$ be a rate-$|\boxn|$ Poisson process which is independent of this family. For each fixed configuration $\sigma\in\Omega^{\boxn}$ and boundary condition $\xi$, define the Markov chain $\tilde{X}_k^{\xi,\sigma}$ so that $\tilde{X}_0^{\xi,\sigma}=\sigma$, while for $k\ge1$ we have $\tilde{X}_k^{\xi,\sigma}(i)=\tilde{X}_{k-1}^{\xi,\sigma}(i)$ for all $i\neq V_k$, and $\tilde{X}_k^{\xi,\sigma}(V_k)=+1$ if $U_{k} \leq p_{V_k}^{\xi}(\tilde{X}_{k-1}^{\xi,\sigma})$ and $\tilde{X}_k^{\xi,\sigma}(V_k)=-1$ if $U_{k} > p_{V_k}^{\xi}(\tilde{X}_{k-1}^{\xi,\sigma})$. The continuous-time heat-bath process on $\boxn$ with boundary condition $\xi$ and initial state $\sigma$ is the process $(X_t^{\xi,\sigma})_{t\geq 0}$ defined by $X_t^{\xi,\sigma}:=\tilde{X}_{N_t}^{\xi,\sigma}$.
\end{definition}

We note that for any fixed boundary condition $\xi$, the family $(X_t^{\xi,\sigma})_{t\ge0,\sigma\in\Omega^{\boxn}}$ defines a grand coupling which is monotone with respect to the natural partial order on $\Omega^{\boxn}$; that is to say, if $\sigma,\tau\in\Omega^{\boxn}$ satisfy $\sigma\le \tau$, then $X_t^{\xi,\sigma}\le X_t^{\xi,\tau}$ for all $t\ge0$. This monotonicity will be very useful in what follows. 

We next introduce the generator and semigroup associated with the heat-bath process. 
For $x \in \boxn$ let $\sigma^{x} \in \Omega^{\boxn}$ denote the unique configuration defined from $\sigma \in \Omega^{\boxn}$ by letting $\sigma^{x}(x) = -\sigma(x)$ and $\sigma^{x}(y) = \sigma(y)$ for all $y\neq x$, i.e. flipping the spin at site $x$. Let $g, h :\Omega^{\boxn} \to \RR$ be arbitrary. We let $\nabla_{x}g(\sigma):= g(\sigma^{x}) - g(\sigma)$ denote the discrete derivative of $g$ at site $x$. The generator, $\igxi = \igl$, for the heat-bath dynamics is then given by
\begin{equation}\label{Eq: generator}
    \ig^{\xi} g(\sigma) = \sum_{x\in \boxn}c_{x}^{\xi}(\sigma)\nabla_{x}g(\sigma)
\end{equation}
where 
\begin{equation}
  c_{x}^{\xi}(\sigma) \coloneqq \frac{1}{2}\left[1-\sigma(x)\tanh\left(\beta S_{x}^{\xi}(\s)\right)\right].
\end{equation}
Since the heat-bath process is reversible, i.e. $\rho_{\boxn, \beta}^{\xi}(\sigma)c_{x}^{\xi}(\sigma) = \rho_{\boxn, \beta}^{\xi}(\sigma^{x})c_{x}^{\xi}(\sigma^{x})$ for all $x \in \boxn$, the generator $\igxi$ is self-adjoint on $L^{2}(\rho_{\boxn, \beta}^{\xi})$, namely
\begin{equation}\label{eq: self adjoint}
    \avgl{g \igxi h} = \avgl{h\igxi g}, \quad \forall\, g, h: \Omega^{\boxn} \to \RR.
\end{equation}
The corresponding semigroup operator $P_t^{\xi} = P_{t, \boxn, \beta}^{\xi}$ is defined by
\begin{equation}
    P_t^{\xi} g(\sigma) = \EE[g(X_{t}^{\xi,\sigma})]. 
\end{equation}
Hille--Yosida theory for Markov semigroups (see, e.g., \cite[Theorem I.2.9]{Liggett2005}) guarantees that the semigroup satisfies the forward and backward Kolmogorov equations
\begin{equation}\label{eq: kolmogorov eqn}
    \frac{\diff}{\dt}P_t^{\xi} g(\sigma) = P_t^{\xi} \igxi g(\sigma) =  \igxi P_t^{\xi} g(\sigma).
\end{equation}
Additionally, the \textit{carr\'e du champ} operator $\gamxi = \gaml$ is the bilinear form satisfying
\begin{equation}\label{Eq: carre du champ}
    \begin{split}
        \gamxi(g, h)(\s) &:= \frac{1}{2}\left(\igxi(gh)(\s) - g(\sigma)\igxi h(\s) - h(\s)\igxi g(\s)\right),\\
        &=\frac{1}{2}\sum_{x \in \boxn}c_{x}^{\xi}(\sigma)\nabla_{x}g(\sigma)\nabla_{x}h(\sigma).
    \end{split}
\end{equation}
By combining~\eqref{eq: self adjoint} and the first equation in~\eqref{Eq: carre du champ} we immediately derive the integration by parts formula
\begin{equation}\label{eq: IBP}
    -\avgl{g\igxi h} = \avgl{\gamxi(g, h)}, \quad \forall\, g, h: \Omega^{\boxn} \to \RR,
\end{equation}
because $\avgl{\igxi g} = 0$ for all $g: \Omega^{\boxn} \to \RR$. The carr\'e du champ operator is central in the theory of Dirichlet forms, and is used to express an $L^{2}$ contraction of the semigroup $P_t^{\xi}$ via the Poincar\'e inequality
\begin{equation}\label{eq: poincare}
    \var_{\boxn, \beta}^{\xi}(g) \leq C_{\beta, d}^{\mathrm{P}}\avgl{\gamxi(g, g)}, \ \ \ \ \forall\, g: \Omega^{\boxn} \to \RR,
\end{equation}
where $C_{\beta, d}^{\mathrm{P}}$ bounds the inverse of the spectral gap for the operator $\igxi$, thus controlling the contraction rate to equilibrium. The logarithmic Sobolev inequality established in~\cite{BauerschmidtDagallier2024} implies the Poincar\'e inequality \eqref{eq: poincare} holds in the entire high-temperature regime with the constant $C_{\beta, d}^{\mathrm{P}}$ being independent of $n$ and of $\xi$.

\subsection{Outline of the paper}
In Section~\ref{sec: steins method}, we outline aspects of Stein's method and introduce the notion of an approximate Stein kernel which we then construct. In Lemma~\ref{lemma: Berry Esseen variance bound} we show the variance of the approximate Stein kernel controls the distance between the laws of the normalised magnetisation $W_{\Lambda}$ and the target Gaussian $Y_{\boxn, \beta}^{\xi}(\Lambda)$. Section~\ref{sec: proof of theorem} proves Theorem~\ref{Theorem: main result 1} assuming the bounds on discrete derivatives of the semigroup in Proposition~\ref{prop: semi group estimates}, which we then establish in Section~\ref{sec: proof of prop}. Section \ref{Sec: stationary gaussian} proves Corollaries~\ref{corol: infinite volume limit theorem} and~\ref{corol: stationary gaussian} using estimates on the difference between the finite-volume and infinite-volume susceptibilities, which we establish in Section~\ref{sec: sus bounds} using standard techniques.

\section{Stein's method for the Kolmogorov distance}\label{sec: steins method} 
Let us begin this section by introducing Stein's method for normal approximation. Recall~\cite{Olver2020} that if a function $f:\RR\to\RR$ is piecewise $\mathrm{C}^1$, then its left derivative exists everywhere. For such a function, $f$, we will adopt the convention throughout that $f'$ denotes its left derivative. Define
\begin{equation}
   \sC:=\{f:\RR\to\RR \mid f \text{ is continuous, piecewise $\mathrm{C}^1$, and } \|f'\|_\infty<\infty\}.
\end{equation}
For $\Lambda \subseteq \boxn$, let $\mathcal{A}_{\Lambda}^{\xi} = \mathcal{A}_{\Lambda, \boxn, \beta}^{\xi}$ be the differential operator such that for all $f\in\sC$, 
\begin{equation}
    \mathcal{A}_{\Lambda}^{\xi}f(x) = \chi_{\boxn, \beta}^{\xi}(\Lambda)f'(x) - xf(x).
\end{equation}
The pair $(\sC, \mathcal{A}_{\Lambda}^{\xi})$ is a characterising pair for $Y_{\boxn, \beta}^{\xi}(\Lambda) \sim \mathcal{N}(0, \chi_{\boxn, \beta}^{\xi}(\Lambda))$ in the sense that
\begin{equation}\label{eq: characterisation}
    \EE\left[\mathcal{A}_{\Lambda}^{\xi}f\left(Y\right)\right] = 0 \ \ \ \forall f\in \sC \iff Y\sim \mathcal{N}\left(0, \chi_{\boxn, \beta}^{\xi}(\Lambda)\right).
\end{equation}
The operator $\mathcal{A}_{\Lambda}^{\xi}$ is related to the generator of an Ornstein-Uhlenbeck process
with invariant measure $\mathcal{N}(0, \chi_{\boxn, \beta}^{\xi}(\Lambda))$
by $\ig^{\mathrm{OU}}f(x) = \mathcal{A}_{\Lambda}^{\xi}f'(x)$ for all suitable $f$.

We then solve the Stein equation for our target distribution
\begin{equation}\label{Eq: Stein Equation}
    \mathcal{A}_{\Lambda}^{\xi} f  = \ind(\cdot \leq z) - \PP(Y_{\boxn, \beta}^{\xi}(\Lambda)\leq z)
\end{equation}
which yields a solution $f = f_{z} \in \sC$ (see, e.g., \cite{ChenShao2005}) given by
\begin{equation}
\begin{split}
    f_{z}(x) &\coloneqq \frac{1}{\chi_{\boxn,\beta}^{\xi}(\Lambda)\,p(x)}\int_{-\infty}^{x}\left(\ind(t \le z) - \PP(Y_{\boxn, \beta}^{\xi}(\Lambda)\leq z) \right)\,p(t)\, \dt,\\
    &\phantom{:}=-\frac{1}{\chi_{\boxn,\beta}^{\xi}(\Lambda)\,p(x)}\int_{x}^{\infty}\left(\ind(t \le z) - \PP(Y_{\boxn, \beta}^{\xi}(\Lambda)\leq z) \right)\,p(t)\, \dt,
\end{split}
\end{equation}
where $p$ is the density of $Y_{\boxn, \beta}^{\xi}(\Lambda)$.

By substituting the rescaled magnetisation $W_{\Lambda}(\sigma)$ into~\eqref{Eq: Stein Equation} and taking expectations (with respect to $\gibbsl$) we find
\begin{equation}\label{Eq: stein equality}
\begin{split}
    \distk(W_{\Lambda}, Y_{\boxn, \beta}^{\xi}(\Lambda))&:=\sup_{z \in \RR}\left|\gibbsl(W_{\Lambda} \leq z) - \PP(Y_{\boxn, \beta}^{\xi}(\Lambda)\leq z)\right|,\\ &= \sup_{z\in\RR}\left|\left\langle \mathcal{A}_{\Lambda}^{\xi} f_{z}(W_{\Lambda})\right\rangle_{\boxn, \beta}^{\xi
    }\right|.
\end{split}
\end{equation}
The Kolmogorov distance, $\distk(W_{\Lambda}, Y_{\boxn, \beta}^{\xi}(\Lambda))$, is then estimated by controlling the last expression in \eqref{Eq: stein equality}, whose dependence on $Y_{\boxn, \beta}^{\xi}(\Lambda)$ is now implicit via the choice of $\mathcal{A}_{\Lambda}^{\xi}$ and the corresponding solution $f_{z}$. To control the last term in~\eqref{Eq: stein equality}, we will need estimates on the solution of the Stein equation. The following result is standard.
\begin{proposition}\label{prop: bounds on stein solution}
    For all $\beta\in[0,\beta_c)$, all nonempty $\Lambda\subseteq\boxn$, every boundary condition $\xi$, and all $u, x, z \in \RR$
    \begin{equation}\label{eq: bound f'}
        |f_{z}'(x)|\leq 1/\chi_{\boxn, \beta}^{\xi}(\Lambda),
    \end{equation}
    \begin{equation}\label{eq: bound xf(x)}
        |xf_{z}(x)| \leq 1,
    \end{equation}
    and 
    \begin{equation}\label{eq: bound f}
        |(x+u)f_z(x+u) - x f_{z}(x)| \leq \frac{|u|}{\chi_{\boxn, \beta}^{\xi}(\Lambda)}\left(|x| + \sqrt{\chi_{\boxn, \beta}^{\xi}(\Lambda)}\right).
    \end{equation}
\end{proposition}
\begin{proof}
    If $\beta = 0$, then $\chi_{\boxn, \beta}^{\xi}(\Lambda) = 1$ for all $\Lambda\neq \emptyset$, and so $Y_{\boxn, 0}^{\xi}(\Lambda)$ is a standard normal. This case can be found in \cite[Lemma 2.2]{ChenShao2005} using the fact that $\sqrt{2\pi}/4 \leq 1$. The general case follows since $f_{z}(x) = (1/\sigma_{\beta})f_{z/\sigma_{\beta}, 0}\left(x/\sigma_{\beta}\right)$ where $\sigma_{\beta} = \sqrt{\chi_{\boxn, \beta}^{\xi}(\Lambda)}$ and $f_{z, 0}$ is the solution to the $\beta = 0$ case.
\end{proof}
A Stein kernel (or Stein coefficient) for a univariate, mean-zero random variable $X$ over the test class $\sC$ is any random variable $\tau_{X}$, defined on the same probability space, such that
\begin{equation}
    \EE[Xf(X)] = \EE[\tau_X f'(X)], \quad \quad \forall f\in \sC.
\end{equation}
The existence of Stein kernels for univariate random variables with suitably well-behaved densities was established in~\cite{Stein1986}. If in addition $X$ has finite variance and $\tau_X \approx \EE[\tau_X] = \var(X)$, one should expect that $X$ is approximately Gaussian since
\begin{equation}
    \EE[Xf(X)] \approx \var(X)\EE[ f'(X)], \quad \forall f\in \sC,
\end{equation}
implying the analogue of~\eqref{eq: stein identity} approximately holds. Quantifying the variance of $\tau_X$ will then quantify the distance between the law of $X$ and a corresponding normal distribution. 

The approach to establishing Theorem~\ref{Theorem: main result 1} is to construct an \textit{approximate Stein kernel}, $\tau_{\Lambda}$, for $W_{\Lambda}$. We may then quantify the distance to a corresponding normal random variable by controlling both the variance of $\tau_{\Lambda}$ \textit{and} the error incurred from the approximation. The form of the approximate Stein kernel is based on a construction of true Stein kernels in \cite{Courtrade2019} which applied to random vectors with densities whose law satisfies a Poincar\'e inequality. Before we define $\tau_{\Lambda}$, we state a useful proposition that will be used in our construction, and which may be thought of as a substitute for the Poincar\'e inequality needed in \cite{Courtrade2019}.
\begin{proposition}\label{prop: decay of semigroup}
    Let $0\leq \beta < \beta_{c}$. There exist constants $C_{\beta,d}>1$ and $c_{\beta, d}\in(0,1]$, independent of $n$, such that for all $t>0$
    \begin{equation}\label{eq: exp contraction}
    \begin{split}
          \sup_{\sigma, x, \xi}\left|\EE[X_t^{\xi,\sigma}(x)] - \avgl{\sigma(x)}\right| \leq C_{\beta,d}\,e^{-c_{\beta, d}t},
    \end{split}
  \end{equation}
  where the supremum is over all $\xi\in\{-1,0,1\}^{\ZZ^d}\cup\{\mathrm{per}\}$, $\sigma\in\Omega^{\boxn}$ and $x\in\boxn$.
\end{proposition}
\begin{remark}\label{rmk: boundary}
  For periodic boundary conditions, Equation~\eqref{eq: exp contraction} was established in~\cite[Lemma 2.1]{Lubetzky2016}. For the remaining boundary conditions, the result follows by combining the uniform-in-volume log-Sobolev inequality~\cite[Theorem 1.1]{BauerschmidtDagallier2024} with~\cite[Theorem 3.1]{Yoshida1997}; the latter being an extension of~\cite[Theorem 1.8]{StroockZegarlinski1992}.
  The argument is essentially immediate when $\xi\in\{-1,1\}^{\ZZ^d}$. 
  Some additional work is required to treat the general case of $\xi\in \{-1, 0, 1\}^{\ZZ^{d}}$, however; the details are provided in Appendix~\ref{sec: decay of semigroup}. 
\end{remark}

As an immediate consequence of Proposition~\ref{prop: decay of semigroup} we have for all initial configurations $\sigma \in \Omega^{\boxn}$
\begin{equation}\label{eq: bound on PtW}
\begin{split}
    |P_t^{\xi}W_{\Lambda}(\sigma)| &\leq\frac{1}{|\Lambda|^{1/2}}\sum_{z\in \Lambda}\left|\EE[X_t^{\xi,\sigma}(z)] - \avgl{\sigma(z)}\right| \\
    & \leq C_{\beta,d}|\Lambda|^{1/2}\,e^{-c_{\beta, d}t}.
\end{split}
\end{equation}
Now define the function $U_{\Lambda} : \Omega^{\boxn} \to \RR$ via
\begin{equation}\label{eq: U def}
    U_{\Lambda}(\sigma) = \int_{0}^{\infty}P_t^{\xi}W_{\Lambda}(\sigma)\, \dt, 
\end{equation}
which is well defined by~\eqref{eq: bound on PtW}. It is easily verified using~\eqref{eq: kolmogorov eqn} that $U_{\Lambda}$ is a solution to the Poisson equation
\begin{equation}\label{eq: discrete stein equation}
    -\igxi U_{\Lambda}(\sigma) = W_{\Lambda}(\sigma),
\end{equation}
since
\begin{equation}
  \begin{split}
    -\igxi U_{\Lambda}(\sigma)
    &= -\sum_{x\in \boxn}c_{x}^{\xi}(\sigma)\left(\int_{0}^{\infty}P_t^{\xi}W_{\Lambda}(\sigma^{x})\, \dt
      -\int_{0}^{\infty}P_t^{\xi}W_{\Lambda}(\sigma)\, \dt\right),\\
    &= -\int_{0}^{\infty}\sum_{x\in \boxn}c_{x}^{\xi}(\sigma)\,\nabla_{x} P_t^{\xi} W_{\Lambda}(\s)\, \dt,\\
    &= -\int_{0}^{\infty}\igxi P_t^{\xi} W_{\Lambda}(\sigma)\, \dt,\\
    &=-\int_{0}^{\infty}\frac{\diff}{\dt}P_t^{\xi} W_{\Lambda}(\sigma)\, \dt,\\
    &=W_{\Lambda}(\sigma) - \lim_{T \to \infty}P_T^{\xi}W_{\Lambda}(\sigma)=W_{\Lambda}(\sigma).
  \end{split}
\end{equation}
where the last line follows from~\eqref{eq: bound on PtW}. 

We may then define the random variable $\tau_{\Lambda}$ via
\begin{equation}
    \tau_{\Lambda}(\sigma):= \gamxi(U_{\Lambda}, W_{\Lambda})(\sigma).
\end{equation}
We will refer to $\tau_{\Lambda}$ as an approximate Stein kernel for $W_{\Lambda}$. 
\begin{lemma}\label{lemma: stein kernel}
    For all $\beta\in[0,\beta_{c})$, all nonempty $\Lambda\subseteq\boxn$, every boundary condition $\xi$, and all $f\in\sC$ we have
    \begin{equation}\label{eq: approx stein kernel}
        \avgl{W_{\Lambda} f(W_{\Lambda})} = \avgl{\tau_{\Lambda} f'(W_{\Lambda})} + \avgl{R(f)}
    \end{equation}
    where
    \begin{equation}\label{eq: Rf def}
        R(f) \coloneq \frac{1}{2}\sum_{x\in \boxn}c_{x}^{\xi}\,\nabla_{x}U_{\Lambda}\,\nabla_{x} W_{\Lambda}\int_{0}^{1}\left(f'(W_{\Lambda} + t\nabla_{x}W_{\Lambda}) - f'(W_{\Lambda})\right)\, \dt.
    \end{equation}
    Moreover,
    \begin{equation}\label{eq: expected value of tau}
        \avgl{\tau_{\Lambda}} = \avgl{W_{\Lambda}^2} = \chi_{\boxn, \beta}^{\xi}(\Lambda).
    \end{equation}
\end{lemma}
\begin{proof}
    Let $f\in\sC$. By~\eqref{eq: discrete stein equation} and the integration by parts formula~\eqref{eq: IBP} we have
    \begin{equation}\label{eq: stein kernel}
        \begin{split}
            \avgl{f(W_{\Lambda}) W_{\Lambda}} & = -\avgl{f(W_{\Lambda}) \igxi U_{\Lambda}} = \avgl{\gamxi(U_{\Lambda}, f(W_{\Lambda}))}.
        \end{split}
    \end{equation}
    Now, since every $f\in\sC$ is absolutely continuous, the fundamental theorem of calculus implies 
    \begin{equation}\label{eq: FTC f'}
    \begin{split}
        \delx (f\circ W_{\Lambda}) &= f(W_{\Lambda} + \delx W_{\Lambda}) - f(W_{\Lambda}),\\
        &=\delx W_{\Lambda} f'(W_{\Lambda}) +\delx W_{\Lambda}\int_{0}^{1}\left(f'(W_{\Lambda} + t\nabla_{x}W_{\Lambda}) - f'(W_{\Lambda})\right)\, \dt.
    \end{split}
    \end{equation}
    And by the definition of $\gamxi(U_{\Lambda}, f(W_{\Lambda}))$, it follows from \eqref{eq: FTC f'} that
    \begin{equation}
        \begin{split}
            \gamxi(U_{\Lambda}, f(W_{\Lambda})) &= \frac{1}{2}\sum_{x \in \boxn}c_{x}^{\xi}\,\delx U_{\Lambda}\,\delx (f\circ W_{\Lambda})\\
            &= \frac{1}{2}\sum_{x \in \boxn}c_{x}^{\xi}\,\delx U_{\Lambda}\,\delx W_{\Lambda}\, f'(W_{\Lambda}) + R(f),\\
            &= \gamxi(U_{\Lambda}, W_{\Lambda})f'(W_{\Lambda}) + R(f),\\
            &=\tau_{\Lambda} f'(W_{\Lambda})+R(f).
        \end{split}
    \end{equation}
    Substituting the last equation into \eqref{eq: stein kernel} yields \eqref{eq: approx stein kernel}.

    To establish~\eqref{eq: expected value of tau}, we substitute $f(x) = x$ into~\eqref{eq: approx stein kernel} and notice that $R(f) = 0$ in this case.
\end{proof}

\begin{remark}\label{remark: stein coupling}
    The equation~\eqref{eq: stein kernel} can be thought of as a continuous-time analogue of a Stein coupling which was introduced in \cite{chen2010steincouplingsnormalapproximation}. If $W_{\Lambda}' = W_{\Lambda}(\sigma')$ and $\sigma'$ is constructed from $\sigma$ by taking a step in the discrete-time heat-bath chain and $U_{\Lambda}(\sigma) = \sum_{n\geq 0}P^{n}W_{\Lambda}(\sigma)$, where $P$ is the transition matrix, this would construct a Stein coupling $(G, W_{\Lambda}, W_{\Lambda}')$ with $G = (U_{\Lambda}'-U_{\Lambda})/2$, where $U_{\Lambda}'$ is defined analogously to $W_{\Lambda}'$. We find it convenient to work in continuous time later in our analysis but additionally, the method of constructing $\tau_{\Lambda}$ in continuous time would also apply to fields with continuous spins by instead using the Langevin dynamics, and $\tau_{\Lambda}$ in that case would indeed be a proper Stein coefficient/kernel analogous to those derived in \cite{Courtrade2019}. 
\end{remark}
Using this formulation, we are able to derive explicit Berry--Esseen bounds for $W_{\Lambda}$, in terms of the variance of $\tau_{\Lambda} = \gamxi(U_{\Lambda}, W_{\Lambda})$. Namely we establish the following estimate.

\begin{lemma}\label{lemma: Berry Esseen variance bound}
    For all $\beta\in[0,\beta_{c})$, all nonempty $\Lambda\subseteq\boxn$ and every boundary condition $\xi$ we have
    \begin{equation}\label{eq: kolmogorov bound var tau}
        \begin{split}
            & \sup_{z\in \RR}\left|\rho_{\boxn, \beta}^{\xi}(W_{\Lambda} \leq z) - \PP\left(Y_{\boxn, \beta}^{\xi}(\Lambda) \leq z\right)\right|\\
            & \quad \quad \quad \quad\quad\quad \leq A_{1}(\Lambda)\sqrt{\var_{\boxn, \beta}^{\xi}(\tau_{\Lambda})} + A_{2}(\Lambda),
        \end{split}
    \end{equation}
    where
    \begin{equation}
        A_{1}(\Lambda) = \frac{1}{\chi_{\boxn, \beta}^{\xi}(\Lambda)}\left(1+\frac{2}{\sqrt{|\Lambda|\chi_{\boxn, \beta}^{\xi}(\Lambda)}}\right),
    \end{equation}
    and
    \begin{equation}
        A_{2}(\Lambda) = \frac{6}{\sqrt{|\Lambda|\chi_{\boxn, \beta}^{\xi}(\Lambda)}}.
    \end{equation}
\end{lemma}
\begin{proof}
    Let $\Phi(z) = \PP\left(Y_{\boxn, \beta}^{\xi}(\Lambda) \leq z\right)$. We fix $0\leq \beta<\beta_{c}$, $n\geq 3$, $\Lambda \subseteq \boxn$ and the boundary condition $\xi$ then omit $\beta$, $\boxn$, $\Lambda$ and $\xi$ from our notation. 
    Fix $x \in \boxn$ and $\sigma \in \Omega^{\boxn}$. For brevity, we let $X_t$ and $X_t^x$ denote the heat-bath processes defined by Definition~\ref{def: Glauber dynamics} with initial states $\sigma$ and $\sigma^x$.
    \begin{equation}
        \nabla_{x}P_tW(\sigma) = \EE[W(X_t^x) - W(X_t)].
    \end{equation}
    By monotonicity of the coupling, if $\sigma^{x}\geq \sigma$ then $X_{t}^{x}\geq X_{t}$ for all $t \geq 0$. Suppose $x\in \Lambda$. Since $W$ is an increasing function, $\delx W(\sigma)\geq 0$ implies $\nabla_{x}P_tW(\sigma)\geq 0$. Similarly, $\delx  W(\sigma)\leq 0$ implies $\nabla_{x}P_tW(\sigma)\leq 0$. On the other hand if $x \in \boxn\setminus \Lambda$ then $\delx W(\sigma) = 0$. It follows that for all $x \in \boxn$ and $\sigma \in \Omega^{\boxn}$ that $\delx W(\sigma)\delx P_tW(\sigma)\geq 0$, and so in particular
    \begin{equation}
        \delx U(\sigma)\delx W(\sigma)\geq 0.
    \end{equation}

    Now, for each $x\in \boxn$ we define
    \begin{equation}\label{eq: KU non negative}
        K_{U}^{x}(s) := \delx U \left(\ind(0 < s \leq \delx W) - \ind(\delx  W \leq s < 0) \right)\geq 0,
    \end{equation}
    which equals $0$ if $x \in \boxn \setminus \Lambda$. Since $|\nabla_{x}W| \leq \delta :=2|\Lambda|^{-1/2}$, $K_{U}^{x}(s)$ satisfies
    \begin{equation}\label{eq: int KU}
        \int_{-\infty}^{\infty}K_{U}^{x}(s)\, \ds = \int_{-\delta}^{\delta}K_{U}^{x}(s)\, \ds = \delx U\, \delx W,
    \end{equation}
    Let $f_{z}$ be the solution to the Stein equation~\eqref{Eq: Stein Equation} corresponding to $z \in \RR$. Since $f_{z}\in\sC$, it follows from~\eqref{eq: approx stein kernel} that, for all $z\in\RR$,
    \begin{equation}
        \begin{split}
            \rho(W\leq z) - \PP(Y\leq z) &= \avg{\mathcal{A} f_{z}(W)}=\avg{(\chi - \tau)f_z'(W)} - \avg{R(f_z)},
        \end{split}
    \end{equation}
    where $R(f)$ is defined in~\eqref{eq: Rf def}. Recalling that $\chi = \avg{\tau}$, Equation~\eqref{eq: bound f'} implies
    \begin{equation}\label{eq: W Y bound}
        |\rho(W\leq z)-\PP(Y\leq z)| \leq \frac{1}{\chi}\sqrt{\var(\tau)} + |\avg{R(f_z)}|.
    \end{equation}
    
    We now establish bounds for $|\avg{R(f_z)}|$. Using~\eqref{eq: KU non negative} and the fact that $f_z$ solves \eqref{Eq: Stein Equation}, we rewrite $R(f_z)$ as
    \begin{equation}\label{eq: R identity R1R2}
        \begin{split}
            R(f_{z}) &= \frac{1}{2}\sum_{x\in \boxn}c_{x}^{\xi}\int_{-\delta}^{\delta}\left[f_{z}'(W+s) - f_{z}'(W)\right]K_{U}^{x}(s)\, \ds,\\
            &=\frac{1}{2\chi}\sum_{x\in \boxn}c_{x}^{\xi}\int_{-\delta}^{\delta}\left[(W+s)f_{z}(W+s) - Wf_{z}(W)\right]K_{U}^{x}(s)\, \ds\\
            &\quad\quad+\frac{1}{2\chi}\sum_{x\in \boxn}c_{x}^{\xi}\int_{-\delta}^{\delta}\left[\ind(W+s\leq z) - \ind(W\leq z)\right]K_{U}^{x}(s)\, \ds,\\
            &=: R_{1}(z) + R_{2}(z).
        \end{split}
    \end{equation}
    
    To bound $|\avg{R_{1}(z)}|$ we may apply~\eqref{eq: bound f} with $|s|\leq \delta$, as well as \eqref{eq: int KU} to obtain
    \begin{equation}\label{eq: R1z bound}
    \begin{split}
        |\avg{R_{1}(z)}| &\leq\avg{\frac{1}{2\chi}\sum_{x\in \boxn}c_{x}^{\xi}\int_{-\delta}^{\delta}\left|(W+s)f_{z}(W+s) - Wf_{z}(W)\right|K_{U}^{x}(s)\, \ds}\\
        &\leq \frac{\delta}{\chi^{2}}\avg{(|W|+\sqrt{\chi})\tau},\\
        &= \frac{\delta}{\chi^2}\avg{|W|(\tau - \chi)} + \frac{\delta \avg{|W|}}{\chi} + \frac{\delta}{\sqrt{\chi}},\\
        &\leq \frac{\delta}{\chi^{3/2}}\sqrt{\var(\tau)} + \frac{2\delta}{\sqrt{\chi}}, 
    \end{split}
    \end{equation}
    where in the last inequality we have applied the Cauchy-Schwarz and Jensen inequalities. For $|\avg{R_{2}(z)}|$, since $\ind(\cdot \leq z)$ is decreasing
    \begin{equation}\label{eq: upper lower R2}
        -\frac{1}{\chi}\avg{\tau \ind(z-\delta<W\leq z)} \leq \avg{R_{2}(z)}\leq \frac{1}{\chi}\avg{\tau \ind(z<W\leq z+\delta)}.
    \end{equation}
    Now define
    \begin{equation}
        I_{z}(x) := \begin{cases}
            0,& x\leq z,\\
            x-z, & z<x\leq z+\delta,\\
            \delta , & x>z+\delta,
        \end{cases}
    \end{equation}
    and note that $I_z\in\sC$. Also note that the bounds $|I_{z}(x)|\leq \delta$ and $I_{z}'(x) = \ind(z<x\leq z+\delta)$ hold for all $x\in\RR$.
    Therefore, by~\eqref{eq: approx stein kernel}, we have
    \begin{equation}\label{eq: tau equality}
        \begin{split}
            \frac{1}{\chi}&\avg{\tau \ind(z<W\leq z+\delta)}\\&=\frac{1}{\chi}\avg{W I_{z}(W)}\\
                          & \quad \quad
                            -\avg{\frac{1}{2\chi}\sum_{x\in \boxn}c_x^{\xi}
                            \int_{-\delta}^{\delta}[\ind(z< W+s \leq z+\delta) - \ind(z<W\leq z+\delta)]K_{U}^{x}(s)\, \ds}\\
                          &=\frac{1}{\chi}\avg{W I_{z}(W)} + \avg{R_{2}(z)} - \avg{R_{2}(z+\delta)}.
        \end{split}
    \end{equation}
    The upper bound in~\eqref{eq: upper lower R2} as well as~\eqref{eq: tau equality} together imply
    \begin{equation}\label{eq: R2 upper unif}
        \avg{R_{2}(z+\delta)} \leq \frac{1}{\chi}\avg{WI_{z}(W)} \leq \frac{\delta}{\sqrt{\chi}}.
    \end{equation}
    An identical argument applied to the lower bound in~\eqref{eq: upper lower R2} yields
    \begin{equation}\label{eq: R2 lower unif}
        \avg{R_{2}(z-\delta)} \geq -\frac{1}{\chi}\avg{WI_{z-\delta}(W)} \geq -\frac{\delta}{\sqrt{\chi}}.
    \end{equation}
    Since the bounds are uniform in $z \in \RR$, it follows that~\eqref{eq: R2 upper unif} and~\eqref{eq: R2 lower unif} imply
    \begin{equation}\label{eq: R2z bound}
        \sup_{z\in \RR}|\avg{R_{2}(z)}| \leq \frac{\delta}{\sqrt{\chi}}.
    \end{equation}
    Finally, combining~\eqref{eq: W Y bound}, \eqref{eq: R identity R1R2}, \eqref{eq: R1z bound} and \eqref{eq: R2z bound} gives
    \begin{equation}
        |\rho(W\leq z)-\PP(Y\leq z)| \leq \frac{1}{\chi}\left(1+\frac{\delta}{\sqrt{\chi}}\right)\sqrt{\var(\tau)} + \frac{3\delta}{\sqrt{\chi}},
    \end{equation}
    and~\eqref{eq: kolmogorov bound var tau} readily follows, recalling that $\delta = 2|\Lambda|^{-1/2}$. 
\end{proof}

The method for obtaining upper bounds in Theorem~\ref{Theorem: main result 1} will then proceed via Lemma~\ref{lemma: Berry Esseen variance bound}, where it will suffice to control $\var_{\boxn, \beta}^{\xi}(\tau_{\Lambda})$. In order to obtain matching lower bounds we will require the following proposition. Its proof is elementary.
\begin{proposition}\label{prop: lattice correction}
    Let $\alpha, \underline{\chi}, \overline{\chi}>0$. Let $Z\sim \mathcal{N}(0, \chi)$ with $0 < \underline{\chi}\leq \chi\leq \overline{\chi}<\infty$. For any $m \in \ZZ^{+}$ and $\gamma \in \RR$ let $X$ be a random variable whose distribution is supported on the translated lattice $\frac{\alpha}{\sqrt{m}}(\ZZ + \gamma)$. Then there exists a positive constant $c = c(\alpha, \underline{\chi}, \overline{\chi})$ such that
    \begin{equation}
        \distk(X, Z) = \sup_{z \in \RR}|\PP(X\leq z) - \PP(Z\leq z)| \geq \frac{c}{\sqrt{m}}.
    \end{equation}
\end{proposition}
\begin{proof}
    Without loss of generality, assume $\gamma\in (0, 1]$. Let $a>0$, $\mathbb{L}:=a(\ZZ + \gamma)$ and define
    \begin{equation}
        I:=\begin{cases}
            (-a/2, 0] &  \text{if }\gamma\in (0, 1/2]\\
            (0, a/2] &  \text{if }\gamma\in (1/2, 1].
        \end{cases}
    \end{equation}
    Then $I\cap \mathbb{L} = \emptyset$, and so if $X$ is an $\mathbb{L}$-valued random variable $\PP(X\in I) = 0$. Hence, setting $I = (u, v]$, if $Z\sim \mathcal{N}(0, \chi)$ we have
    \begin{equation}
        \begin{split}
            \PP(Z\in I) &= |\PP(Z\in I) - \PP(X\in I)|,\\
            &\leq |\PP(Z\leq v) - \PP(X\leq v)|+|\PP(Z\leq u) - \PP(X\leq u)|,\\
            &\leq 2\distk(X, Z)
        \end{split}
    \end{equation}
    By symmetry, $\PP(Z\in I) = \PP(0<Z\leq a/2)$. Setting $a=\alpha/\sqrt{m}$, $m\in \ZZ^{+}$ and using $\underline{\chi}\leq \chi\leq \overline{\chi}$ implies
    \begin{equation}
        \begin{split}
            \distk(X, Z)&\geq \frac{1}{2}\PP(Z\in I)\\
            &=\frac{1}{2}\frac{1}{\sqrt{2\pi \chi}}\int_{0}^{\alpha/2\sqrt{m}}\exp\left(-\frac{x^{2}}{2\chi}\right)\, \dx,\\
            &\geq\frac{\alpha}{4\sqrt{m}}\frac{1}{\sqrt{2\pi \chi}}\,\exp\left(-\frac{\alpha^{2}}{8m\chi}\right),\\
            &\geq \frac{c}{\sqrt{m}},
        \end{split}
    \end{equation}
    where
    \begin{equation}\label{eq: lower bound complete theorem}
        c = \frac{\alpha}{4\sqrt{2\pi \overline \chi}}\,\exp\left(-\frac{\alpha^{2}}{8\underline{\chi}}\right).
    \end{equation}
\end{proof}

\section{Proof of Theorem \ref{Theorem: main result 1}}\label{sec: proof of theorem}
In order to prove the upper bound in Theorem~\ref{Theorem: main result 1} it suffices, given Lemma~\ref{lemma: Berry Esseen variance bound}, to control $\var_{\boxn, \beta}^{\xi}(\tau_{\Lambda})$. Our strategy for doing so is to apply a Poincar\'e inequality, and then use estimates on discrete derivatives of $U_{\Lambda}$ defined in~\eqref{eq: U def}. The estimates on $U_{\Lambda}$ follow immediately from the following bounds on the discrete derivatives of the semigroup $P_t^{\xi} W_{\Lambda}$. In what follows we let $\dist(x, y)$ denote the graphical distance between the vertices $x$ and $y$ on the graph $\BB_{n}^{d}$ (or $\TT_{n}^{d}$ if $\xi = \mathrm{per}$).
\begin{proposition}\label{prop: semi group estimates}
    For all $\beta\in[0,\beta_{c})$, boundary conditions $\xi$, initial states $\sigma \in \Omega^{\boxn}$, nonempty $\Lambda\subseteq\boxn$, sites $x, y \in \boxn$ and $t\geq 0$, the following hold:
    \begin{equation}\label{eq: 1st deriv PtW bound}
        |\delx P_t^{\xi}W_{\Lambda}(\sigma)| \leq \frac{C_{1}}{\sqrt{|\Lambda|}}\exp(-c_{1} t),
    \end{equation}
    \begin{equation}\label{eq: 2nd deriv PtW bound}
        |\dely \delx P_t^{\xi}W_{\Lambda}(\sigma)| \leq \frac{C_{2}}{\sqrt{|\Lambda|}}\exp(-c_{2} t)\exp(-k_2\, \dist(x, y)),
    \end{equation}
    where $c_{1}, c_{2}, C_{1}, C_{2}, k_2 \in (0, \infty)$ depend only on $\beta$ and $d$.
\end{proposition}

We delay the proof of Proposition~\ref{prop: semi group estimates} to Section~\ref{sec: proof of prop} since this is the main technical result. The desired bounds on the discrete derivatives of $U_{\Lambda}$ now follow as an immediate corollary of Proposition~\ref{prop: semi group estimates}.
\begin{corollary}\label{cor: disc deriv U}
  For all $\beta\in[0,\beta_{c})$, boundary conditions $\xi$, initial states $\sigma \in \Omega^{\boxn}$, nonempty $\Lambda\subseteq\boxn$, and $x, y \in \boxn$, the following hold:
    \begin{equation}\label{eq: delx U}
        |\nabla_{x}U_{\Lambda}(\sigma)| \leq \frac{u_1}{\sqrt{|\Lambda|}},
    \end{equation}
    \begin{equation}\label{eq: dely delx U}
        |\nabla_{y}\nabla_{x}U_{\Lambda}(\sigma)| \leq \frac{u_2}{\sqrt{|\Lambda|}}\exp(-k_{2}\,\dist(x, y)),
    \end{equation}
    where $u_1 := C_1/c_1$ and $u_2 := C_2/c_2$, with $c_{1}, c_2, C_1, C_{2}$ and $k_{2}$ defined in Proposition~\ref{prop: semi group estimates}.
\end{corollary}
\begin{proof}
    The estimates \eqref{eq: delx U} and \eqref{eq: dely delx U} are straightforward. For \eqref{eq: delx U} we use \eqref{eq: 1st deriv PtW bound} to find
    \begin{equation}
        \begin{split}
            |\delx U_{\Lambda}(\sigma)| &\leq \int_{0}^{\infty}\left|\delx P_t^{\xi}W_{\Lambda}(\sigma)\right|\, \dt,\\
            & \leq \frac{C_{1}}{c_{1}\sqrt{|\Lambda|}},
        \end{split}
    \end{equation}
    and \eqref{eq: dely delx U} is obtained from \eqref{eq: 2nd deriv PtW bound} similarly.
\end{proof}

In addition, we will require the following bounds on the finite-volume susceptibility.
\begin{proposition}\label{prop: susceptibility bounds}
    For all $\beta\in[0,\beta_{c})$, all nonempty $\Lambda\subseteq \boxn$, and all boundary conditions $\xi$,
    \begin{equation}\label{eq: susceptibility bounds}
        \frac{1}{\cosh(2d\beta)^2}\leq \chi_{\boxn, \beta}^{\xi}(\Lambda) \leq \chi_{\ZZ^{d}, \beta}.
    \end{equation}
    Moreover, if $\xi \in \{0, \mathrm{per}\}$, then $\chi_{\boxn, \beta}^{\xi}(\Lambda)\geq 1$.
\end{proposition}
\begin{proof}
    Recall that
    \begin{equation}
        \chi_{\boxn, \beta}^{\xi}(\Lambda) = \frac{1}{|\Lambda|}\sum_{x, y \in \Lambda}\left(\avgl{\sigma(x) \sigma(y)} - \avgl{\sigma(x)}\avgl{\sigma(y)}\right).
    \end{equation}
    Since $\avgl{\sigma(x) \sigma(y)} - \avgl{\sigma(x)}\avgl{\sigma(y)} \geq 0$ by the FKG inequality \cite{FKG1971}, we have
    \begin{equation}\label{eq: lower bound chi}
        \chi_{\boxn, \beta}^{\xi}(\Lambda) \geq \frac{1}{|\Lambda|}\sum_{x \in \Lambda}\left(1 - \left(\avgl{\sigma(x)}\right)^{2}\right).
    \end{equation}
    Now, it is easily verified that $\lvert\avgl{\sigma(x)|(\sigma(y))_{y\neq x}}\rvert = \lvert\tanh(\beta S_{x}^{\xi}(\sigma))\rvert \leq \lvert\tanh(2d\beta)\rvert$. Applied to~\eqref{eq: lower bound chi} this implies
    \begin{equation}
        \begin{split}
            \chi_{\boxn, \beta}^{\xi}(\Lambda) &\geq\frac{1}{|\Lambda|}\sum_{x \in \Lambda}\left(1 - \left(\avgl{\avgl{\sigma(x)|\{\sigma(y)\}_{y\neq x}}}\right)^{2}\right)\\
            &\geq1-\tanh(2d\beta)^{2}\\
            &= \frac{1}{\cosh(2d\beta)^2}
        \end{split}
    \end{equation}
    which establishes the lower bound in~\eqref{eq: susceptibility bounds}. In the case when $\xi \in \{0, \mathrm{per}\}$ we observe that $\avgl{\sigma(x)} = 0$ for all $x \in \Lambda$. Thus~\eqref{eq: lower bound chi} immediately implies $\chi_{\boxn, \beta}^{0}(\Lambda) \geq 1$ and $\chi_{\boxn, \beta}^{\mathrm{per}}(\Lambda) \geq 1$ for every $\Lambda \subseteq \boxn$.

    For the upper bound, first consider periodic boundary conditions $\xi = \mathrm{per}$. By \cite[Theorem 2.1]{Liu2025}
    \begin{equation}
        \avg{\sigma(0)\sigma(x)}_{\boxn, \beta}^{\mathrm{per}} \leq \sum_{u\in \ZZ^{d}}\avg{\sigma(0)\sigma(x+nu)}_{\ZZ^{d}, \beta}.
    \end{equation}
    Since $\avg{\sigma(x)}_{\boxn, \beta}^{\mathrm{per}} = 0$ for all $x\in \boxn$ and $\avg{\sigma(x)\sigma(y)}_{\boxn, \beta}^{\mathrm{per}} \geq 0$ for all $x, y \in \boxn$, by translation invariance under periodic boundary conditions
    \begin{equation}\label{eq: sus periodic bound}
        \begin{split}
            \chi_{\boxn, \beta}^{\mathrm{per}}(\Lambda) &= \frac{1}{|\Lambda|}\sum_{x, y \in \Lambda}\avg{\sigma(x)\sigma(y)}_{\boxn, \beta}^{\mathrm{per}}\\
            &\leq \frac{1}{|\Lambda|}\sum_{\substack{x \in \boxn\\y\in \Lambda}}\avg{\sigma(x)\sigma(y)}_{\boxn, \beta}^{\mathrm{per}}\\
            &=\sum_{x\in \boxn}\avg{\sigma(0)\sigma(x)}_{\boxn, \beta}^{\mathrm{per}}\\
            & \leq \sum_{x\in \boxn}\sum_{u \in \ZZ^{d}}\avg{\sigma(0)\sigma(x + nu)}_{\ZZ^{d}, \beta}\\
            &=\chi_{\ZZ^{d}, \beta}.
        \end{split}
    \end{equation}
    Now consider $\xi \neq \mathrm{per}$. We have that
    \begin{equation}\label{eq: xi bound by per}
    \begin{split}
        \avgl{\sigma(x) \sigma(y)} - \avgl{\sigma(x)}\avgl{\sigma(y)} &\leq \avg{\sigma(x) \sigma(y)}_{\boxn, \beta}^{0} \leq \avg{\sigma(x) \sigma(y)}_{\boxn, \beta}^{\mathrm{per}},
    \end{split}
    \end{equation}
    where the first inequality follows by \cite[Corollary 1.3]{Ding2023} and the last inequality follows by Griffiths' inequality \cite[Theorem 2]{griffiths1967correlations}. By summing over $x, y \in \Lambda$ and dividing by $|\Lambda|$, Equations~\eqref{eq: xi bound by per} and~\eqref{eq: sus periodic bound} then together imply $\chi_{\boxn, \beta}^{\xi}(\Lambda) \leq \chi_{\boxn, \beta}^{\mathrm{per}}(\Lambda) \leq \chi_{\ZZ^{d}, \beta}$. 
\end{proof}

By combining Lemma~\ref{lemma: Berry Esseen variance bound}, Corollary~\ref{cor: disc deriv U} and Proposition~\ref{prop: susceptibility bounds}, we may now prove Theorem~\ref{Theorem: main result 1}.

\begin{proof}[Proof of Theorem \ref{Theorem: main result 1}]
    Let us begin with the upper bound in~\eqref{eq: berry esseen same variance}. We estimate the variance of $\tau_{\Lambda} = \gamxi(U_{\Lambda}, W_{\Lambda})$ by applying the Poincar\'e inequality~\eqref{eq: poincare} to $\tau_{\Lambda}$. Namely
    \begin{equation}\label{eq: poincare eq}
        \var_{\boxn, \beta}^{\xi}(\tau_{\Lambda}) \leq C_{\beta, d}^{\mathrm{P}}\avgl{\gamxi(\tau_{\Lambda}, \tau_{\Lambda})},
    \end{equation}
    with constant $C_{\beta, d}^{\mathrm{P}}$ independent of $n$ and $\xi$.  We begin by computing $(\dely \tau_{\Lambda})^2$. Fix a configuration $\sigma$, and note that, for $y\in \boxn$,
    \begin{equation}\label{eq: sum decomposition}
      \begin{split}
        \dely \tau_{\Lambda}(\sigma)
        &= \frac{1}{2}\sum_{x\in \boxn}\dely\left[c_{x}^{\xi}(\sigma)\,\delx W_{\Lambda}(\sigma)\,\delx U_{\Lambda}(\sigma) \right],\\
        &=\frac{1}{2}\sum_{x\in \boxn}\dely c_{x}^{\xi}(\sigma)\, \delx W_{\Lambda}(\sigma)\,\delx U_{\Lambda}(\sigma)\\
        & \qquad + \frac{1}{2}\sum_{x\in \boxn} c_{x}^{\xi}(\sy)\, \dely\delx W_{\Lambda}(\sigma)\,\delx U_{\Lambda}(\sigma)\\
        & \qquad + \frac{1}{2}\sum_{x\in \boxn} c_{x}^{\xi}(\sy)\, \delx W_{\Lambda}(\sigma^{y})\,\dely\delx U_{\Lambda}(\sigma),\\
        &=: J_{1}(y)+ J_2(y) + J_3(y).
        \end{split}
    \end{equation}
     Direct calculation shows
    \begin{equation}
         \delx W_{\Lambda}(\sigma) = -\frac{2\sigma(x)}{|\Lambda|^{1/2}}\ind(x\in \Lambda),
     \end{equation}
     \begin{equation}\label{eq: dely delx W}
         \dely \delx W_{\Lambda}(\sigma) = \frac{4\sigma(x)}{|\Lambda|^{1/2}}\ind(y =x)\ind(x\in \Lambda),
     \end{equation}
     and
     \begin{equation}\label{eq: dely cx}
     \begin{split}
         \dely c_{x}^{\xi}(\sigma) &= \sigma(x)\tanh(\beta S_{x}^{\xi}(\sigma))\ind(y=x)\\
         & \quad\quad\quad \quad+\frac{\sigma(x)}{2}\left[\tanh(\beta S_{x}^{\xi}(\sigma)) - \tanh(\beta S_{x}^{\xi}(\sigma) -2\beta\sigma(y))\right]\ind(y \sim x).
     \end{split}
     \end{equation}
     We may then replace the sums over $x\in\boxn$ by sums over $x\in \Lambda$ in~\eqref{eq: sum decomposition}. Since $\gamxi(\tau_{\Lambda}, \tau_{\Lambda}) = \frac{1}{2}\sum_{y\in \boxn}c_{y}^{\xi}\,(J_{1}(y)+J_{2}(y)+J_{3}(y))^2$, we proceed by bounding the squares of the $J_{i}$.

     For $J_{1}$, we first observe that by~\eqref{eq: dely cx} 
     \begin{equation}\label{eq: dely cx bound}
         |\dely c_x^{\xi}| \leq \ind(\dist(x, y)\leq 1), \quad\quad x, y \in \boxn.
     \end{equation}
     We may apply $|\delx W_{\Lambda}| \leq 2|\Lambda|^{-1/2}$ and Equations~\eqref{eq: delx U} and~\eqref{eq: dely cx bound} to obtain
     \begin{equation}\label{eq: S1}
         J_{1}(y)^2 \leq \left(\frac{1}{2}\sum_{x\in \Lambda: \dist(x, y)\leq 1} \, \frac{2u_{1}}{|\Lambda|}\right)^{2} \leq \frac{(2d+1)^2 u_{1}^{2}}{|\Lambda|^2}\ind(\dist(y, \Lambda)\leq 1),
     \end{equation}
     where in the last inequality we have used the fact that $x$ has at most $2d$ neighbours and that the sum is empty if $\dist(y, \Lambda)>1$. For $J_{2}$ we apply~\eqref{eq: dely delx W}, \eqref{eq: delx U} and $|c_{x}^{\xi}| \leq 1$ to get
     \begin{equation}\label{eq: S2}
         J_{2}(y)^2 \leq \left(\frac{1}{2} \, \frac{4u_{1}}{|\Lambda|}\ind(y\in \Lambda)\right)^{2} = \frac{4u_{1}^2}{|\Lambda|^2}\ind(y\in \Lambda).
     \end{equation}
     Finally for $J_{3}$ we apply~\eqref{eq: dely delx U}
     \begin{equation}\label{eq: S3}
     \begin{split}
         J_{3}(y)^2 &\leq \left(\frac{1}{2}\sum_{x\in \Lambda}\frac{2u_2}{|\Lambda|}e^{-k_{2}\,\dist(x, y)}\right)^{2},\\
         &=\frac{u_2^{2}}{|\Lambda|^2}\sum_{x_{1}, x_{2}\in \Lambda}e^{-k_2(\dist(x_{1}, y)+\dist(x_{2}, y))},\\
         &\leq \frac{u_2^{2}}{|\Lambda|^2}\sum_{\substack{x_{1}\in \Lambda\\x_2\in \ZZ^{d}}}e^{-k_2(\dist(x_{1}, y)+|y-x_{2}|_{1})},\\
         &=\frac{u_2^{2}\coth(k_{2}/2)^{d}}{|\Lambda|^2}\sum_{x_{1}\in \Lambda}e^{-k_2\,\dist(x_{1}, y)},
     \end{split}
     \end{equation}
    where for the second inequality we have used that $\dist(x_{2}, y) = |x_{2}' - y|_{1}$ for some $x_{2}'\in \ZZ^{d}$, for either the torus graphical distance or the box graphical distance (in which case $x_{2}' = x_{2}$). Now since
    \begin{equation}
        \sum_{y\in \boxn}\ind(\dist(y, \Lambda)\leq 1)\leq (2d+1)|\Lambda|,
    \end{equation}
    and
    \begin{equation}
        \sum_{\substack{y\in \boxn\\x_{1}\in \Lambda}}e^{-k_{2}\dist(x_{1}, y)}\leq\sum_{\substack{y\in \ZZ^{d}\\x_{1}\in \Lambda}}e^{-k_{2}|y-x_{1}|_{1}} = |\Lambda|\coth(k_{2}/2)^d,
    \end{equation}
    combining~\eqref{eq: S1}, \eqref{eq: S2} and \eqref{eq: S3} with $|c_{y}^{\xi}|\leq 1$ yields
     \begin{equation}\label{eq: gamb tau bound}
         \begin{split}
             \gamxi(\tau_{\Lambda}, \tau_{\Lambda}) &=\frac{1}{2}\sum_{y\in \boxn}c_{y}^{\xi}\,[J_{1}(y)+J_{2}(y)+J_{3}(y)]^2,\\
             & \leq \frac{3}{2}\sum_{y\in \boxn}[J_{1}(y)^2 + J_{2}(y)^2 + J_{3}(y)^2],\\
             &\leq \frac{K_{\beta, d}}{|\Lambda|},
         \end{split}
     \end{equation}
     where
     \begin{equation}
         K_{\beta, d} =\frac{3}{2}\left((2d+1)^3u_1^2+4u_1^2+u_2^2\coth(k_2/2)^{2d}\right).
     \end{equation}
     
     Applying \eqref{eq: poincare eq} we find that $\sqrt{\var_{\boxn, \beta}^{\xi}(\tau_{\Lambda})} \leq \sqrt{C_{\beta, d}^{\mathrm{P}}K_{\beta, d}/|\Lambda|}$. Combined with Lemma~\ref{lemma: Berry Esseen variance bound} we obtain
     \begin{equation}\label{eq: dkol prelim bound}
     \begin{split}
       \sup_{z\in \RR}\left|\rho_{\boxn, \beta}^{\xi}(W_{\Lambda} \leq z) - \PP(Y_{\boxn, \beta}^{\xi}(\Lambda)\leq z)\right|
       &\leq A_{1}(\Lambda)\sqrt{C_{\beta, d}^{\mathrm{P}}K_{\beta, d}/|\Lambda|}+ A_{2}(\Lambda). 
     \end{split}
     \end{equation}
     By Proposition~\ref{prop: susceptibility bounds} we can find constants $a_1$ and $a_2$ which depend only on $\beta$ and $d$ such that
     \begin{equation}\label{eq: a1 bound chi}
         A_{1}(\Lambda) = \frac{1}{\chi_{\boxn, \beta}^{\xi}(\Lambda)}\left(1+\frac{2}{\sqrt{|\Lambda|\chi_{\boxn, \beta}^{\xi}(\Lambda)}}\right) \leq a_1,
     \end{equation}
     and
     \begin{equation}\label{eq: a2 bound}
         A_{2}(\Lambda) = \frac{6}{\sqrt{|\Lambda|\chi_{\boxn, \beta}^{\xi}}} \leq \frac{a_2}{\sqrt{|\Lambda|}}.
     \end{equation}
     Combining~\eqref{eq: dkol prelim bound}, \eqref{eq: a1 bound chi} and \eqref{eq: a2 bound} we can find a constant $C = C(\beta, d)$ such that
     \begin{equation}
         \distk(W_{\Lambda}, Y_{\boxn, \beta}^{\xi}(\Lambda)) = \sup_{z\in \RR}\left|\rho_{\boxn, \beta}^{\xi}(W_{\Lambda} \leq z) - \PP(Y_{\boxn, \beta}^{\xi}(\Lambda)\leq z)\right| \leq \frac{C}{\sqrt{|\Lambda|}}.
     \end{equation}

     For the lower bound let $\underline{\chi} = 1/\cosh(2d\beta)^2$ and $\overline{\chi} = \chi_{\ZZ^{d}, \beta}$. Then for any $n\geq 3$, Proposition~\ref{prop: susceptibility bounds} implies $\chi_{\boxn, \beta}^{\xi}(\Lambda) \in [\underline{\chi}, \overline{\chi}]$ and moreover, $W_{\Lambda}$ is supported on $|\Lambda|^{-1/2}(\ZZ - \avgl{S_{\Lambda}(\sigma)})$. Applying Proposition~\ref{prop: lattice correction} then yields
     \begin{equation}
         \distk(W_{\Lambda}, Y_{\boxn, \beta}^{\xi}(\Lambda)) = \sup_{z\in \RR}\left|\rho_{\boxn, \beta}^{\xi}(W_{\Lambda} \leq z) - \PP(Y_{\boxn, \beta}^{\xi}(\Lambda)\leq z)\right|\geq \frac{c}{\sqrt{|\Lambda|}},
     \end{equation}
     where
     \begin{equation}\label{eq: cdef}
         c = c(\beta, d) = \frac{1}{4\sqrt{2\pi \chi_{\ZZ^{d}, \beta}}}\exp(-\cosh(2d\beta)^{2}/8).
     \end{equation}
     This completes the proof.
\end{proof}

\section{Proof of Proposition \ref{prop: semi group estimates}}\label{sec: proof of prop}
For given $\sigma\in\Omega^{\boxn}$ and $x\in\boxn$, let $(X_t^{\xi,\sigma,x})_{t\ge0}$ denote the process defined in Definition~\ref{def: Glauber dynamics} with boundary condition $\xi$ and initial state $\sigma^x$.
When $\xi$ and $\sigma$ are clear from context we will at times omit them from our notation.

We now summarise some key terminology and definitions introduced in \cite{Lubetzky2013, Lubetzky2016}.
Let $0<T_1<T_2<\ldots$ denote the jump times of the Poisson process $(N_t)_{t\ge0}$.
The \emph{update sequence} from $0$ to $t$, corresponding to the grand coupling given by Definition~\ref{def: Glauber dynamics}, is the set of tuples $(V_i, U_i, T_i)$ with $T_{i}\leq t$. For two sites $u, v \in \boxn$ we say that $u$ \textit{reaches} $v$ via the update sequence up to time $t$, and write $u\overset{\xi, t}{\rightsquigarrow}v$, if $u=v$ or there exists integer $m\in \ZZ^{+}$ and a subset of the update sequence $\{(V_{i_{j}}, U_{i_j}, T_{i_{j}})\}_{1\leq j\leq m}$ such that $u \sim V_{i_{1}}\sim \cdots\sim V_{i_{m}} = v$ and $0\leq T_{i_{1}} \leq \cdots \leq T_{i_{m}} \leq t$. We emphasise that the adjacency relation differs depending on whether or not $\xi=\mathrm{per}$. 
For two subsets $A, B \subseteq \boxn$, we then say $A$ reaches $B$ via the update sequence up to time $t$, denoted by $A\overset{\xi, t}{\rightsquigarrow}B$, if there exists $u\in A$ and $v\in B$ such that $u\overset{\xi, t}{\rightsquigarrow}v$. 

For $A \subseteq \boxn$, write
\begin{equation}
X_t^{\xi,\sigma}(A) \coloneqq \bigl(X_t^{\xi,\sigma}(x)\bigr)_{x\in A}.
\end{equation}
This leads naturally to the definition of the \textit{update functions} which are analogous to the quantities defined in \cite{Lubetzky2016}. We adopt a notation which is slightly more convenient for our purposes.

\begin{definition}[Forward Update Function]\label{eq: forward update}
    For a given set $A\subseteq \boxn$ the forward update function is the random set, $\mathcal{F}^{\xi}(A, t)$, that consists of all sites $u\in \boxn$ such that $A\overset{\xi, t}{\rightsquigarrow}u$. 
\end{definition}

\begin{definition}[Backward Update Function]\label{eq: backward update}
    For a given set $A\subseteq \boxn$ the backward update function is the random set, $\mathcal{B}^{\xi}(A, t)$, that consists of all sites $u\in \boxn$ such that $u\overset{\xi, t}{\rightsquigarrow}A$.
\end{definition}

\begin{definition}[Update Support Function]\label{def: update support}
    Given a boundary condition $\xi$, the update sequence, and a set $A\subseteq \boxn$, we define the update support function via
    \begin{equation*}
    \mathcal{S}^{\xi}(A, t):=\{x\in \boxn: \exists\, \sigma\in \Omega^{\boxn} \text{ such that } X_t^{\xi,\sigma}(A)\neq X_{t}^{\xi,\sigma,x}(A)\}.
    \end{equation*} 
\end{definition}

We write $\mathcal{S}^{\xi}(v, t) = \mathcal{S}^{\xi}(\{v\}, t)$ for $v\in \boxn$, with analogous conventions for the forward and backward update functions. Note that $\mathcal{S}^{\xi}(v, t) \subseteq \mathcal{B}^{\xi}(v, t)$ for all $v \in \boxn$ and $t\geq 0$, and all boundary conditions $\xi$.

We will require the following lemma concerning the temporal and spatial decay of the update functions, which was proved in \cite[Lemma 2.1]{Lubetzky2016} for periodic boundary conditions using \cite[Theorem 3.1]{Martinelli1994a} and the finite speed of information flow. The proof for the general case follows essentially verbatim by applying Proposition~\ref{prop: decay of semigroup}. Let $\dist(x, y)$ denote the graphical distance between the vertices $x$ and $y$ corresponding to the graph $\BB_{n}^{d}$ (or $\TT_{n}^{d}$ if $\xi = \mathrm{per}$), and define
\begin{equation}
    B_{r}(v):=\{x\in \boxn: \dist(v, x)<r\},\quad \quad v\in \boxn.
\end{equation}

\begin{lemma}\label{lemma: sup and upd}
  For all $\beta\in[0,\beta_{c})$, all boundary conditions $\xi$, and all $t\ge0$, the update functions satisfy the following. For all $v\in \boxn$ we have
    \begin{equation}\label{eq: Fsup bound}
        \PP(\mathcal{S}^{\xi}(v, t) \neq \emptyset) \leq C_{\beta,d}\,e^{-c_{\beta, d}t},
    \end{equation}
    where $c_{\beta, d}, C_{\beta,d}$ are as given in Proposition~\ref{prop: decay of semigroup}. Moreover, for any $r>20dt$,
    \begin{equation}\label{eq: Fupd bound}
        \PP(\mathcal{F}^{\xi}(v, t) \not\subset B_{r}(v)) \leq e^{-r}.
    \end{equation}
    \begin{equation}\label{eq: Bupd bound}
        \PP(\mathcal{B}^{\xi}(v, t) \not\subset B_{r}(v)) \leq e^{-r}.
    \end{equation}
\end{lemma}

\begin{proof}
  Fix a boundary condition $\xi$. The stated results are immediate when $t=0$, so suppose instead that $t>0$. 

  We begin by proving~\eqref{eq: Fsup bound}.
  Let $X_t^{\xi,+}$ denote the heat-bath process whose initial state has every spin set to $+1$, and
  let $X_t^{\xi,-}$ denote the process whose initial state has every spin set to $-1$.
  Due to the monotonicity of the grand coupling, given an update sequence, $\mathcal{S}^{\xi}(v, t) \neq \emptyset$ if and only if
$X_{t}^{\xi,+}(v) \neq X_t^{\xi,-}(v)$. Hence
    \begin{equation}
    \begin{split}
      \PP(\mathcal{S}^{\xi}(v, t) \neq \emptyset)
      &= \PP[X^{\xi,+}_{t}(v)\neq X_{t}^{\xi,-}(v)],\\
        &=\frac{1}{2}\EE[X^{\xi,+}_{t}(v) - X_{t}^{\xi,-}(v)],\\
      &  = \frac{1}{2}\EE\left[X^{\xi,+}_{t}(v)- \avgl{\sigma(v)}\right]
        - \frac{1}{2}\EE\left[X_{t}^{\xi,-}(v) - \avgl{\sigma(v)}\right]\\
        & \leq C_{\beta,d} \,e^{-c_{\beta, d}t}
    \end{split}
    \end{equation}
    where the final inequality follows by Proposition~\ref{prop: decay of semigroup}.

    We next consider~\eqref{eq: Fupd bound}. Let
    \begin{equation}
        \mathcal{W}_{r}(v) \coloneqq \{w = (w_{0}, w_1, \dots, w_{\lceil r\rceil}): w_{0} = v \text{ and } \dist(w_i, w_{i-1}) = 1\}
    \end{equation}
    be the set of walks of length $k\coloneqq\lceil r\rceil$ started at $v$. For $\mathcal{F}^{\xi}(v, t)\not \subset B_{k}(v)$, there must be a vertex $u\in B_{k}(v)^{c}$ such that $v \overset{\xi, t}{\rightsquigarrow}u$. For this to occur, there must be a sequence of times $0 < t_{1}\leq \cdots\leq t_{k} \leq t$ and $w\in \mathcal{W}_{r}(v)$ such that $w_{i}$ is to be updated at time $t_{i}$. For fixed $w\in \mathcal{W}_{r}$ call this event $A_{r}(w)$.
Conditional on \(N_t\), for \(A_r(w)\) to occur there must exist indices $1\le i_1<\cdots<i_k\le N_t$
such that \(V_{i_j}=w_j\) for every \(1\le j\le k\). The union bound then implies
\begin{equation}\label{eq: cond Ar}
\mathbb P(A_r(w)\mid N_t) \le \binom{N_t}{k}|\Lambda_n|^{-k},
\end{equation}
and so
\begin{equation}\label{eq: Ar}
\mathbb P(A_r(w))
\le
|\Lambda_n|^{-k}\mathbb E\binom{N_t}{k}
=
\frac{t^k}{k!}
\le
\left(\frac{et}{k}\right)^k
=
e^{-k(\log(k/t)-1)}.
\end{equation}
    Since there are at most $(2d)^{k}$ such walks, applying the union bound once more and using $k=\lceil r\rceil \geq r >20dt$ yields
    \begin{equation}
        \PP(\mathcal{F}^{\xi}(v, t)\not \subset B_{\lceil r\rceil}(v)) \leq (2d)^{k} e^{-k(\log(k/t)- 1)} \leq e^{-r}.
    \end{equation}
    But since the graph distance is integer valued, $B_{\lceil r \rceil}(v)=B_r(v)$, and so~\eqref{eq: Fupd bound} follows.

For~\eqref{eq: Bupd bound} we may argue similarly. For $\mathcal{B}^{\xi}(v, t)\not\subset B_{k}(v)$, there must be a vertex $u\in B_{k}(v)^{c}$ such that $u \overset{\xi, t}{\rightsquigarrow}v$. For this to occur, there must be a sequence of times $0 < t_{1}\leq \cdots\leq t_{k} \leq t$ and $w\in \mathcal{W}_{r}(v)$ such that $w_{k-i}$ is to be updated at time $t_{i}$. Call this event $A^{*}_{r}(w)$. Clearly $A^{*}_{r}(w)$ has the same law as $A_{r}(w)$ and so~\eqref{eq: cond Ar} and~\eqref{eq: Ar} also hold for $A^{*}_{r}(w)$. Equation~\eqref{eq: Bupd bound} then follows from the same argument as above.\qedhere
\end{proof}

For $\sigma\in\Omega^{\boxn}$ and distinct $x,y\in\boxn$, let $\sigma^{x, y}$ denote the configuration obtained from $\sigma$ by flipping the spins at $x$ and $y$, and let $(X_t^{\xi,\sigma,x,y})_{t\ge0}$ denote the process defined in Definition~\ref{def: Glauber dynamics} with boundary condition $\xi$ and initial state $\sigma^{x,y}$.
For $\eta, \omega \in \Omega^{\boxn}$, let
\begin{equation}
  D_{\Lambda}(\eta, \omega) \coloneqq \{z\in \Lambda : \eta(z)\neq \omega(z)\},
\end{equation}
that is, the set of disagreeing sites in $\Lambda$ between configurations $\eta$ and $\omega$. We will also require the following lemma.
\begin{lemma}\label{lemma: disagreements}
  Fix $\beta\in[0,\beta_{c})$, nonempty $\Lambda\subseteq\boxn$, boundary condition $\xi$, initial state $\sigma\in \Omega^{\boxn}$,  and distinct $x, y\in \boxn$. 
For any $t\ge0$, if $D_{\Lambda}(X_{t}^{\xi,\sigma}, X_{t}^{\xi,\sigma,x}) \neq D_{\Lambda}(X_{t}^{\xi,\sigma,y}, X_t^{\xi,\sigma,x, y})$, then there exists $z \in \Lambda$ such that $x, y \in \mathcal{S}^{\xi}(z, t)$.
\end{lemma}
\begin{proof}

  We fix a boundary condition $\xi$ and configuration $\sigma\in\Omega^{\boxn}$, and suppress them from our notation.
  Now fix distinct $x,y\in\boxn$. We prove the contrapositive.
  Assume, therefore, that for all $z\in \Lambda$ we have $\{x, y\}\not \subseteq \mathcal{S}(z, t)$. For every $z \in \Lambda$, exactly one of the following then holds:
    \begin{enumerate}[label=(\roman*)]
        \item $x \notin \mathcal{S}(z, t)$\label{item: x not in},
        \item $x\in \mathcal{S}(z, t)$ and $y \notin \mathcal{S}(z, t)$\label{item: y not in}.
    \end{enumerate}
Suppose~\ref{item: x not in} holds for a given $z\in \Lambda$. Then by definition $X_{t}(z) = X_{t}^{x}(z)$ and $X_{t}^{y}(z) = X_{t}^{x,y}(z)$, since in each case the initial configurations differ only at $x$. This implies $z \notin D_{\Lambda}(X_{t}, X_{t}^{x})$ and $z \notin D_{\Lambda}(X_{t}^{y}, X_{t}^{x, y})$. On the other hand, suppose~\ref{item: y not in} holds for a given $z\in \Lambda$. Then by definition $X_{t}(z) = X_{t}^{y}(z)$ and $X_{t}^{x}(z) = X_{t}^{x, y}(z)$, which implies that $X_{t}(z) - X_{t}^{x}(z) = X_{t}^{y}(z) - X_{t}^{x, y}(z)$. Thus, either $z \notin D_{\Lambda}(X_{t}, X_{t}^{x})$ and $z \notin D_{\Lambda}(X_{t}^{y}, X_{t}^{x, y})$, or $z \in D_{\Lambda}(X_{t}, X_{t}^{x})$ and $z \in D_{\Lambda}(X_{t}^{y}, X_{t}^{x, y})$. We conclude $D_{\Lambda}(X_{t}, X_{t}^{x}) = D_{\Lambda}(X_{t}^{y}, X^{x, y}_{t})$.
\end{proof}

With Lemma~\ref{lemma: sup and upd} and Lemma~\ref{lemma: disagreements}, we may now prove Proposition~\ref{prop: semi group estimates}.

\begin{proof}[Proof of Proposition~\ref{prop: semi group estimates}]
  We begin with the proof of \eqref{eq: 1st deriv PtW bound}. Fix $\sigma \in \Omega^{\boxn}$ and $x \in \boxn$.
  Since \(|\nabla_x h(\sigma^x)|=|\nabla_x h(\sigma)|\) for any \(h:\Omega^{\Lambda_n}\to\mathbb R\), we may assume,
  without loss of generality, that \(\sigma(x)=-1\), so that $\sigma^{x} \geq \sigma$. From the monotonicity of the grand coupling we have $X_{t}^{\xi,\sigma,x}\geq X_{t}^{\xi,\sigma}$ for all $t \geq 0$. We then have for every $t\geq 0$
    \begin{equation}\label{eq: delx PtW rewrite}
        \begin{split}
          \delx P_t^{\xi} W_{\Lambda}(\sigma) &= \frac{1}{\sqrt{|\Lambda|}}\sum_{z\in \Lambda}
                                                \EE[X_{t}^{\xi,\sigma,x}(z) - X_{t}^{\xi,\sigma}(z)],\\
            &=\frac{2}{\sqrt{|\Lambda|}}\sum_{z \in \Lambda}\PP[X_{t}^{\xi,\sigma,x}(z) \neq X_{t}^{\xi,\sigma}(z)],\\
            &\leq\frac{2}{\sqrt{|\Lambda|}}\sum_{z \in \Lambda}\PP[x\in \mathcal{S}^{\xi}(z, t)],
        \end{split}
    \end{equation}
    If $x\in \mathcal{S}^{\xi}(z, t)$, then necessarily $x\in \mathcal{B}^{\xi}(z, t)$ and $\mathcal{S}^{\xi}(z, t)\neq \emptyset$. Combining with~\eqref{eq: delx PtW rewrite}, then applying the Cauchy-Schwarz inequality together with~\eqref{eq: Fsup bound}, we have
    \begin{equation}\label{eq: bound delx PtW upd}
        \begin{split}
            \delx P_t^{\xi}W_{\Lambda}(\sigma)&\leq \frac{2}{\sqrt{|\Lambda|}}\sum_{z \in \Lambda}\PP(\mathcal{S}^{\xi}(z, t)\neq \emptyset, x\in \mathcal{B}^{\xi}(z, t)),\\
            & \leq  \frac{2}{\sqrt{|\Lambda|}}\sum_{z \in \Lambda}\left[\PP(\mathcal{S}^{\xi}(z, t)\neq \emptyset)\PP( x\in \mathcal{B}^{\xi}(z, t))\right]^{1/2}\\
            &\leq \frac{2\sqrt{C_{\beta,d}}\,e^{-c_{\beta, d} t/2}}{\sqrt{|\Lambda|}}\sum_{z \in \Lambda}\PP( x\in \mathcal{B}^{\xi}(z, t))^{1/2}
        \end{split}
    \end{equation}
    We now split the sum of $z$ depending on whether or not $\dist(x, z) > 20dt$. That is
    \begin{equation}
    \begin{split}
      \sum_{z \in \Lambda}\PP( x\in \mathcal{B}^{\xi}(z, t))^{1/2}
      &= \sum_{\substack{z\in \Lambda:\\ \dist(x, z)> 20dt}}\PP( x\in \mathcal{B}^{\xi}(z, t))^{1/2}\\
      & \quad \quad \quad \quad + \sum_{\substack{z\in \Lambda:\\ \dist(x, z)\leq 20dt}}\PP( x\in \mathcal{B}^{\xi}(z, t))^{1/2}\\
      &=: A_1 + A_2.
    \end{split}
    \end{equation}
    For $A_1$, we may apply~\eqref{eq: Bupd bound} to obtain
    \begin{equation}\label{eq: A1 bound coth}
        \begin{split}
            A_{1} &\leq \sum_{\substack{z\in \Lambda:\\ \dist(x, z)> 20dt}} e^{-\dist(x, z)/2}\\
            &\leq \sum_{z \in \ZZ^{d}} e^{-|x-z|_{1}/2}\\
            &= \coth(1/4)^{d},
        \end{split}
    \end{equation}
    where for the second inequality we have used that $\dist(x, z) = |x - z'|_{1}$ for some $z'\in \ZZ^{d}$, for either the torus graphical distance or the box graphical distance (in which case $z' = z$). On the other hand, for $A_{2}$ we simply bound the probability by $1$ to get the number of sites $z\in \Lambda$ with $\dist(x, z)$ at most $20dt$. The number of such sites is bounded by the volume of the continuous $L^1$ ball with radius $ 20dt + d/2$. Thus
    \begin{equation}\label{eq: A2 bound unif}
        A_{2}\leq \frac{2^d}{d!}(20dt + d/2)^{d} \leq e^{d}(40t + 1)^{d},
    \end{equation}
    where in the last line we used the bound $e^d \geq d^d/d!$. Combining~\eqref{eq: bound delx PtW upd}, \eqref{eq: A1 bound coth} and~\eqref{eq: A2 bound unif}, we finally have
    \begin{equation}
        \begin{split}
            |\nabla_{x}P_t^{\xi}W_{\Lambda}(\sigma)| &\leq \frac{2\sqrt{C_{\beta,d}}}{\sqrt{|\Lambda|}}\left(\coth(1/4)^de^{-c_{\beta, d}t/2} + e^{d}(40 t + 1)^de^{-c_{\beta, d}t/2}\right),\\
            & \leq \frac{C_{1}}{\sqrt{|\Lambda|}} \exp(-c_{1}t),
        \end{split}
    \end{equation}
    where
    \begin{equation}
        c_{1} = c_{\beta, d}/4,
    \end{equation}
    and
    \begin{equation}
        C_{1} = 2\sqrt{C_{\beta,d}}\left(\coth(1/4)^{d} + e^{c_{\beta, d}/160}\left(\frac{160 d}{c_{\beta, d}}\right)^{d}\right).
    \end{equation}
    Since $\sigma, \xi$ and $x$ were arbitrary in the above argument, we now conclude that~\eqref{eq: 1st deriv PtW bound} holds for all boundary conditions $\xi$, and all $\sigma \in \Omega^{\boxn}$ and $x\in \boxn$.  
     
    We now turn our attention to~\eqref{eq: 2nd deriv PtW bound}. Fix a boundary condition $\xi$, initial state $\sigma \in \Omega^{\boxn}$, and distinct sites $x, y \in \boxn$, and assume now, without loss of generality, that $\sigma(x) = 1$.
    Due to the monotonicity of the grand coupling we have
    \begin{equation}\label{eq: dely delx Pt equality}
        \begin{split}
          \dely\delx P_t^{\xi} W_{\Lambda}(\sigma)
          &= \frac{1}{\sqrt{|\Lambda|}}\sum_{z \in \Lambda}\EE\left[[X_{t}^{\xi,\sigma}(z) - X_{t}^{\xi,\sigma,x}(z)]-[X_{t}^{\xi,\sigma,y}(z) - X_{t}^{\xi,\sigma,x, y}(z)]\right],\\
          &=\frac{2}{\sqrt{|\Lambda|}}\EE\left[|D_{\Lambda}(X_t^{\xi,\sigma}, X_t^{\xi,\sigma,x})| - |D_{\Lambda}(X_t^{\xi,\sigma,y}, X_t^{\xi,\sigma,x, y})|\right]
        \end{split}
    \end{equation}
    We proceed by first controlling $\EE[|D_{\Lambda}(X_t^{\xi,\sigma}, X_t^{\xi,\sigma,x})|^{2}]$. By the Cauchy-Schwarz inequality
    \begin{equation}
    \begin{split}
      \sqrt{\EE[|D_{\Lambda}(X_t^{\xi,\sigma}, X_t^{\xi,\sigma,x})|^{2}]}
      &= \left(\sum_{z, w \in \Lambda}\PP\left[X_{t}^{\xi,\sigma,x}(z) \neq X_{t}^{\xi,\sigma}(z), X_{t}^{\xi,\sigma,x}(w) \neq X_{t}^{\xi,\sigma}(w)\right]\right)^{1/2},\\
        &\leq \sum_{z\in \Lambda}\PP[X_{t}^{\xi,\sigma,x}(z) \neq X_{t}^{\xi,\sigma}(z)]^{1/2},\\
        & \leq C_{\beta,d}^{1/4}e^{-c_{\beta, d}t/4}\sum_{z\in \Lambda}\PP(x\in \mathcal{B}^{\xi}(z, t))^{1/4}
    \end{split}
    \end{equation}
    where the last line follows from the same argument used in~\eqref{eq: bound delx PtW upd}. By the same argument as above, decomposing into sums over $\dist(x, z) >20dt$ and $\dist(x, z) \leq20dt$ and controlling these sums we have for all $\sigma\in \Omega^{\boxn}$ and all $x\in \boxn$
    \begin{equation}\label{eq: D squared}
        \sqrt{\EE[|D_{\Lambda}(X_t^{\xi,\sigma}, X_t^{\xi,\sigma,x})|^{2}]} \leq M_{1}\exp(-m_1 t)
    \end{equation}
    where
    \begin{equation}
        m_{1} = c_{\beta, d}/8
    \end{equation}
    and 
    \begin{equation}
        M_1 = C_{\beta,d}^{1/4}\left(\coth(1/8)^{d} + e^{c_{\beta, d}/320}\left(\frac{320 d}{c_{\beta, d}}\right)^{d}\right).
    \end{equation}
    
    We now return our attention to~\eqref{eq: dely delx Pt equality}. Fix $t\geq 0$ and let
    \begin{equation}
        A_{t} := D_{\Lambda}(X_t^{\xi,\sigma}, X_t^{\xi,\sigma,x}), \quad B_{t}:=D_{\Lambda}(X_t^{\xi,\sigma,y}, X_t^{\xi,\sigma,x, y}),
    \end{equation} 
    which both satisfy~\eqref{eq: D squared}. Additionally, by Lemma~\ref{lemma: disagreements}, if $A_{t} \neq B_{t}$ there must be a vertex $z\in \Lambda$ such that $x, y \in \mathcal{S}^{\xi}(z, t)$. This implies $x, y \in \mathcal{B}^{\xi}(z, t)$ and hence $z \in \mathcal{F}^{\xi}(x, t)\cap \mathcal{F}^{\xi}(y, t)$. Combined with \eqref{eq: dely delx Pt equality} we can then show using the Cauchy-Schwarz inequality
    \begin{equation}\label{eq: dely delx PtW estimate}
        \begin{split}
          |\dely\delx P_t^{\xi}& W_{\Lambda}(\sigma)|\\
                               &= \left|\frac{2}{\sqrt{|\Lambda|}}\EE\left[\left(|A_{t}| - |B_{t}|\right)\ind(\mathcal{F}^{\xi}(x, t)\cap \mathcal{F}^{\xi}(y, t) \neq \emptyset)\right]\right|,\\
            & \leq \frac{2}{\sqrt{|\Lambda|}}\sqrt{\EE[\left(|A_{t}| - |B_{t}|\right)^{2}]}\, \PP(\mathcal{F}^{\xi}(x, t)\cap \mathcal{F}^{\xi}(y, t) \neq \emptyset)^{1/2}\\
            &\leq \frac{4}{\sqrt{|\Lambda|}}\left(\sqrt{\EE[|A_{t}|^{2}]}+\sqrt{\EE[|B_{t}|^{2}]}\right) \PP(\mathcal{F}^{\xi}(x, t)\cap \mathcal{F}^{\xi}(y, t) \neq \emptyset)^{1/2}\\
            &\leq \frac{8M_{1}}{\sqrt{|\Lambda|}}\exp(-m_1 t)\PP(\mathcal{F}^{\xi}(x, t)\cap \mathcal{F}^{\xi}(y, t) \neq \emptyset)^{1/2},
        \end{split}
    \end{equation}
    where in the last line we applied~\eqref{eq: D squared}. Let $d_{xy}:= \dist(x, y)/2$, and note that if $\mathcal{F}^{\xi}(x, t)\cap \mathcal{F}^{\xi}(y, t) \neq \emptyset$, then either $\mathcal{F}^{\xi}(x, t) \not\subset B_{d_{xy}}(x)$ or $\mathcal{F}^{\xi}(y, t) \not\subset B_{d_{xy}}(y)$. We therefore obtain from the union bound
    \begin{equation}\label{eq: xy Fupd bound}
        \begin{split}
            &\PP(\mathcal{F}^{\xi}(x, t)\cap \mathcal{F}^{\xi}(y, t) \neq \emptyset)\\
            & \quad \quad \quad \quad \leq \PP\left(\mathcal{F}^{\xi}(x, t) \not \subset B_{d_{xy}}(x)\right)+\PP\left(\mathcal{F}^{\xi}(y, t) \not \subset B_{d_{xy}}(y)\right),\\
            &\quad \quad \quad \quad =: P_1 + P_2
        \end{split} 
    \end{equation}
    Combining~\eqref{eq: dely delx PtW estimate} and~\eqref{eq: xy Fupd bound} we arrive at
    \begin{equation}
        |\dely\delx P_t^{\xi} W_{\Lambda}(\sigma)| \leq \frac{8M_{1}}{\sqrt{|\Lambda|}}\exp(-m_1 t)(P_1 + P_2)^{1/2}.
    \end{equation}
    To obtain~\eqref{eq: 2nd deriv PtW bound}, if $d_{xy} =\dist(x, y)/2 > 20dt$ we apply \eqref{eq: Fupd bound} to get
    \begin{equation}
        |\dely\delx P_t^{\xi} W_{\Lambda}(\sigma)| \leq \frac{8\sqrt{2}M_{1}}{\sqrt{|\Lambda|}}\exp(-m_1 t)\exp(-\dist(x, y)/4).
    \end{equation}
    On the other hand if $0<\dist(x, y)/2 \leq 20dt$, then we simply use $P_i\le 1$ to obtain
    \begin{equation}
    \begin{split}
        |\dely\delx P_t^{\xi} W_{\Lambda}(\sigma)|&\leq \frac{8\sqrt{2}M_{1}}{\sqrt{|\Lambda|}}\exp(-m_1 t)\\
        & \leq\frac{8\sqrt{2}M_{1}}{\sqrt{|\Lambda|}}\exp(-m_1 t/2)\exp(-m_{1}\dist(x, y)/80d).
    \end{split}
    \end{equation}
    In either case, we now conclude that for all $x\neq y$
    \begin{equation}
      |\dely\delx P_t^{\xi} W_{\Lambda}(\sigma)| \leq \frac{C_{2}}{\sqrt{|\Lambda|}}\exp(-c_{2} t)\exp(-k_2\dist(x, y)),
      \label{eq: 2nd deriv PtW x!=y}
    \end{equation}
    where
    \begin{equation}
        C_{2} = 2^{7/2}\sqrt{C_{\beta,d}}\,\left(\coth(1/8)^{d} + e^{c_{\beta, d}/320}\left(\frac{320 d}{c_{\beta, d}}\right)^{d}\right),
    \end{equation}
    \begin{equation}
        c_{2} = c_{\beta, d}/16,
    \end{equation}
    and 
    \begin{equation}
        k_{2} = c_{\beta, d}/640d.
      \end{equation}
      Finally, note that $\nabla_x\nabla_x h=-2\nabla_x h$ for any $h:\Omega^{\boxn}\to\RR$. Therefore, in the case $x=y$, it follows as an immediate consequence of~\eqref{eq: 1st deriv PtW bound} that
      \begin{equation}
      |\delx\delx P_t^{\xi} W_{\Lambda}(\sigma)| \leq \frac{2C_{1}}{\sqrt{|\Lambda|}}\exp(-c_{1} t).
      \end{equation}
      But it is easily verified that $2C_1\le C_2$ and $c_2\le c_1$. And so it follows that~\eqref{eq: 2nd deriv PtW x!=y} in fact holds for all $x,y\in\boxn$. And since $\sigma$ and $\xi$ were arbitrary in the above argument, we now conclude that~\eqref{eq: 2nd deriv PtW bound} holds for all boundary conditions $\xi$, all $\sigma \in \Omega^{\boxn}$, and all $x,y\in\boxn$. 
\end{proof}
\section{Proof of Corollaries~\ref{corol: infinite volume limit theorem} and~\ref{corol: stationary gaussian}}\label{Sec: stationary gaussian}
We may immediately establish Corollary~\ref{corol: infinite volume limit theorem} from Theorem~\ref{Theorem: main result 1}.
\begin{proof}[Proof of Corollary~\ref{corol: infinite volume limit theorem}]
    Consider Theorem~\ref{Theorem: main result 1} with free boundary conditions, and fix a finite subset $\Lambda\subset\ZZ^{d}$ and $\beta<\beta_{c}$. Define 
    \begin{equation}
        F_{\boxn, \beta}^{0}(z;\Lambda):=|\rho_{\boxn, \beta}^{0}(W_{\Lambda}\leq z) - \PP(Y_{\boxn, \beta}^{0}(\Lambda)\leq z)|,
    \end{equation}
    and
    \begin{equation}
        F_{\ZZ^{d}, \beta}(z;\Lambda):=|\rho_{\ZZ^{d}, \beta}(W_{\Lambda}\leq z) - \PP(Y_{\ZZ^{d}, \beta}(\Lambda)\leq z)|.
    \end{equation}
    Then by~\eqref{eq: unique Gibbs state} we have 
    \begin{equation}
        \lim_{n \to \infty}F_{\boxn, \beta}^{0}(z;\Lambda) = F_{\ZZ^{d}, \beta}(z;\Lambda).
    \end{equation}
    It follows from Theorem~\ref{Theorem: main result 1} that
    \begin{equation}
        \sup_{z\in \RR}F_{\ZZ^{d}, \beta}(z;\Lambda) \leq \limsup_{n \to \infty}\sup_{z\in \RR}F_{\boxn, \beta}^{0}(z;\Lambda) \leq \frac{C}{\sqrt{|\Lambda|}},
    \end{equation}
    which is precisely the upper bound in~\eqref{eq: berry esseen same variance inf}. The lower bound may be established by applying Proposition~\ref{prop: lattice correction}.
\end{proof}

We now turn our attention to Corollary~\ref{corol: stationary gaussian}. To unify the notation for the susceptibility, for the remainder of this section we define 
\begin{equation}\label{eq: chi infty}
    \chi_{\boxn, \beta}^{\infty}:= \chi_{\ZZ^{d}, \beta}(\boxn) = \frac{1}{n^{d}}\sum_{x, y \in \boxn}\avg{\sigma(x)\sigma(y)}_{\ZZ^{d}, \beta},
\end{equation}
and extend the possible values of $\xi$ to include $\xi = \infty$, representing the infinite volume object. Additionally, we let $\Lambda = \boxn$ and suppress $\Lambda$ in the notation, i.e. we write $W := W_{\boxn} = S_{\boxn}/\sqrt{|\boxn|}$ and $\chi_{\boxn, \beta}^{\xi} = \chi_{\boxn, \beta}^{\xi}(\boxn)$ for $\xi \in \{0, \mathrm{per}\}$. Using Theorem~\ref{Theorem: main result 1} and Corollary~\ref{corol: infinite volume limit theorem}, we can now establish bounds between $W$ and $Y_{\beta}\sim \mathcal{N}(0, \chi_{\ZZ^{d}, \beta})$ for both free and periodic boundary conditions or under the infinite-volume Gibbs measure. We will require the following estimates, which are simple consequences of the Simon-Lieb inequality~\cite{Simon1980, Lieb1980}, the exponential decay of correlations and recently established results in \cite{Liu2025}. 
\begin{lemma}\label{lemma: sus bounds free and periodic}
    For every $0<\beta < \beta_{c}$, there are positive constants $k^{\mathrm{per}}, K^{\mathrm{per}},k^{0},K^{0}, k^{\infty}, K^{\infty}$ (which depend only on $\beta$ and $d$) such that for all $n\geq 3$
    \begin{equation}\label{eq: chi per bound}
         0 \leq \chi_{\ZZ^{d}, \beta} - \chi_{\boxn, \beta}^{\mathrm{per}} \leq K^{\mathrm{per}}e^{-k^{\mathrm{per}}n},
    \end{equation}
    \begin{equation}\label{eq: chi 0 bound}
        \frac{k^{0}}{n} \leq \chi_{\ZZ^{d}, \beta} - \chi_{\boxn, \beta}^{0} \leq \frac{K^{0}}{n},
    \end{equation}
    and 
    \begin{equation}\label{eq: chi inf bound}
        \frac{k^{\infty}}{n} \leq \chi_{\ZZ^{d}, \beta} - \chi_{\boxn, \beta}^{\infty} \leq \frac{K^{\infty}}{n}
    \end{equation}
\end{lemma}
The bounds~\eqref{eq: stationary gaussian per},~\eqref{eq: stationary gaussian free} and~\eqref{eq: stationary gaussian inf} in Corollary~\ref{corol: stationary gaussian} now follow rather easily.

\begin{proof}[Proof of Corollary~\ref{corol: stationary gaussian}]
    Let us begin with the upper bounds. Fix $\xi \in \{\mathrm{per}, 0, \infty\}$. If $\xi \in \{0, \mathrm{per}\}$ let $W=W(\sigma)$ be distributed according to $\rho_{\boxn, \beta}^{\xi}$, and otherwise let $W$ be distributed according to $\rho_{\ZZ^{d}, \beta}$. Let $Y_{\boxn, \beta}^{\xi}\sim \mathcal{N}(0, \chi_{\boxn, \beta}^{\xi})$ and $Y_{\beta} \sim \mathcal{N}(0, \chi_{\ZZ^{d}, \beta})$. Then, by the triangle inequality and Theorem~\ref{Theorem: main result 1} (or Corollary~\ref{corol: infinite volume limit theorem} if $\xi = \infty$), 
    \begin{equation}\label{eq: triangle kol}
        \begin{split}
            \distk(W, Y_{\beta}) &\leq \distk\left(W, Y_{\boxn, \beta}^{\xi}\right) + \distk\left(Y_{\boxn, \beta}^{\xi}, Y_{\beta}\right)\\
            & \leq \frac{C}{n^{d/2}} + \distk\left(Y_{\boxn, \beta}^{\xi}, Y_{\beta}\right).
        \end{split}
    \end{equation}
    An elementary calculation (see, e.g., \cite[Proposition 2.6]{JohnstonProchno2022}) shows that the Kolmogorov distance between two centred Gaussian random variables with variances $\chi_{\ZZ^{d}, \beta}\geq \chi_{\boxn, \beta}^{\xi}\geq 1$ satisfies
    \begin{equation}\label{eq: kol gaussian}
        \distk\left(Y_{\boxn, \beta}^{\xi}, Y_{\beta}\right) \leq \frac{3}{8}\frac{\chi_{\ZZ^{d}, \beta} - \chi_{\boxn, \beta}^{\xi}}{\chi_{\boxn, \beta}^{\xi}} \leq \frac{3}{8}\left(\chi_{\ZZ^{d}, \beta} - \chi_{\boxn, \beta}^{\xi}\right). 
    \end{equation}
    Now combining~\eqref{eq: triangle kol} and~\eqref{eq: kol gaussian} as well as the respective bounds~\eqref{eq: chi per bound},~\eqref{eq: chi 0 bound} and~\eqref{eq: chi inf bound} immediately proves the upper bounds in~\eqref{eq: stationary gaussian per},~\eqref{eq: stationary gaussian free} and~\eqref{eq: stationary gaussian inf} with appropriate constants $C_{\mathrm{per}}, C_{0}, C_{\infty}$. 

    By applying Proposition~\ref{prop: lattice correction} with $\underline{\chi} = \overline{\chi} = \chi_{\ZZ^{d}, \beta}$, if $W$ is distributed according to $\rho_{\boxn, \beta}^{\xi}$ or $\rho_{\ZZ^{d}, \beta}$ then
    \begin{equation}\label{eq: a prioir kol bound}
        \distk(W, Y_{\beta})\geq \frac{c}{n^{d/2}}
    \end{equation}
    with $c = c(\beta, d)$ defined in~\eqref{eq: cdef}. This establishes the lower bound~\eqref{eq: stationary gaussian per} in the case of periodic boundary conditions with $c_{\mathrm{per}} = c$, and the lower bounds in~\eqref{eq: stationary gaussian free} and~\eqref{eq: stationary gaussian inf} when $d=1, 2$ with $c_{0}, c_{\infty} = c$. It remains to show that if $W$ is distributed according to $\rho_{\boxn, \beta}^{0}$ then
    \begin{equation}\label{eq: dgeq3 free}
        \sup_{z\in \RR}|\rho_{\boxn, \beta}^{0}(W\leq z) - \PP(Y_{\beta}\leq z)| \geq \frac{c_{0}}{n}, \quad \quad d\geq 3,
    \end{equation}
    and otherwise if $W$ is distributed according to $\rho_{\ZZ^{d}, \beta}$ then
    \begin{equation}\label{eq: dgeq3 inf}
        \sup_{z\in \RR}|\rho_{\ZZ^{d}, \beta}(W\leq z) - \PP(Y_{\beta}\leq z)| \geq \frac{c_{\infty}}{n}, \quad \quad d\geq 3,
    \end{equation}
    for some positive constants $c_{0} = c_{0}(\beta, d)$ and $c_{\infty} = c_{\infty}(\beta, d)$. We will show~\eqref{eq: dgeq3 free} and note that~\eqref{eq: dgeq3 inf} follows by an identical argument.
    
    First, via an elementary calculation we obtain a lower bound on $\distk\left(Y_{\boxn, \beta}^{0}, Y_{\beta}\right)$. If $Z \sim \mathcal{N}(0, 1)$ and $\alpha:= \sqrt{\chi_{\boxn, \beta}^{0}/\chi_{\ZZ^{d}, \beta}}<1$ then for every $n \geq 3$
    \begin{equation}\label{eq: kol gaussian lower}
        \begin{split}
            \distk\left(Y_{\boxn, \beta}^{0}, Y_{\beta}\right)&\geq \left|\PP\left(Y_{\boxn, \beta}^{0} \leq \sqrt{\chi_{\boxn, \beta}^{0}}\right) - \PP\left(Y_{\beta}\leq \sqrt{\chi_{\boxn, \beta}^{0}}\right)\right|,\\
            &=\PP\left(\alpha < Z \leq 1\right), \\
            &\geq \frac{1}{\sqrt{2\pi}}e^{-1/2}\left(1-\alpha\right),\\
            &\geq \frac{1}{2\sqrt{2\pi}}e^{-1/2}(1-\alpha^{2}),\\
            &= \frac{1}{2\sqrt{2\pi}}e^{-1/2}\frac{1}{\chi_{\ZZ^{d}, \beta}}(\chi_{\ZZ^{d}, \beta} - \chi_{\boxn, \beta}^{0}).
        \end{split}
    \end{equation}
    Combining~\eqref{eq: chi 0 bound} and \eqref{eq: kol gaussian lower} we then have
    \begin{equation}
        \distk\left(Y_{\boxn, \beta}^{0}, Y_{\beta}\right) \geq \frac{c_{0}'}{n}.
    \end{equation}
    Now let $W$ be distributed according to $\rho_{\boxn, \beta}^{0}$. By the reverse triangle inequality and Theorem~\ref{Theorem: main result 1} we have for all $n$ sufficiently large
    \begin{equation}
    \begin{split}
        \distk(W, Y_{\beta}) &\geq \left| \distk\left(Y_{\boxn, \beta}^{0}, Y_{\beta}\right) -\distk\left(W, Y_{\boxn, \beta}^{0}\right)\right|,\\
        & \geq \frac{c_{0}'}{n} -\frac{C}{n^{d/2}}\\
        &\geq \frac{c_{0}'}{2n}
    \end{split}
    \end{equation}
    In light of the prior estimate~\eqref{eq: a prioir kol bound}, for each $d\ge3$ we may then find a constant $c_{0}\leq c_{0}'/2$ sufficiently small that for all $n\geq 3$
    \begin{equation}
        \distk(W, Y_{\beta})\geq \frac{c_{0}}{n}.
    \end{equation}
    This establishes~\eqref{eq: dgeq3 free}, and~\eqref{eq: dgeq3 inf} follows by an identical argument. Combined with~\eqref{eq: a prioir kol bound} when $d=1, 2$, this establishes the lower bounds in~\eqref{eq: stationary gaussian free} and~\eqref{eq: stationary gaussian inf} and completes the proof of the corollary.
\end{proof}

\section{Proof of Lemma~\ref{lemma: sus bounds free and periodic}}\label{sec: sus bounds}
We prove~\eqref{eq: chi per bound},~\eqref{eq: chi 0 bound} and~\eqref{eq: chi inf bound} separately. Before we prove~\eqref{eq: chi per bound}, we introduce some notation and a result from \cite{Liu2025}. Fix the lattice dimension $d\in \ZZ^{+}$ and write $\torusn:=\torusnd = (\boxn, E_{\TT}(\boxn))$. For periodic boundary conditions, we inherit the notion of modular addition which we will define and make use of in this section. Following \cite[Section 2.2]{Liu2025}, given $x, x'\in \ZZ^{d}$ we write $x\cong x'$ if $x-x' \in n\ZZ^{d}$. For any $x\in \ZZ^{d}$ we then write $\pi x$ for the unique point in $\Lambda_{n}$ (corresponding to a vertex in $\torusn$) such that $\pi x \cong x$. Then for every $x, y \in \boxn$, torus addition and subtraction is defined via
\begin{equation}
    x\oplus y := \pi(x+y)\in \boxn, \quad \quad x\ominus y := \pi(x-y)\in \boxn.
\end{equation}

For $x, y \in \boxn $ we let $\tau_{\beta}^{\torusn}(x, y) :=  \avg{\sigma(x)\sigma(y)}_{\boxn, \beta}^{\mathrm{per}}$ be the torus two-point function, and for $x, y \in \ZZ^{d}$ we let $\tau_{\beta}(x, y) = \langle \sigma(x)\sigma(y)\rangle_{\ZZ^{d, }\beta}$ be the infinite-volume two-point function. By translation invariance, $\tau_{\beta}^{\torusn}(x, y) = \tau_{\beta}^{\torusn}(0, y\ominus x)$ and so we write $\tau_{\beta}^{\torusn}(x) := \tau_{\beta}^{\torusn}(0, x)$ (and similarly for $\tau_{\beta}$). Also define the unwrapped two-point function
    \begin{equation}
        T_{\beta}^{(n)}(x) := \sum_{u\in \ZZ^{d}}\tau_{\beta}(x+nu), \quad x\in \boxn,
    \end{equation}
    which satisfies for all $\beta<\beta_{c}$
    \begin{equation}\label{eq: Tbeta sus}
        \sum_{x\in \boxn}T_{\beta}^{(n)}(x) =\sum_{x\in \boxn} \sum_{u\in \ZZ^{d}}\tau_{\beta}(x+nu) = \chi_{\ZZ^{d}, \beta}.
    \end{equation}
    Let $\star$ denote convolution in $\torusn$, i.e. for $x \in \boxn$ we define $(f\star g)(x) := \sum_{y\in \boxn}f(y)g(x\ominus y)$. The following result was proved in \cite[Corollary 2.2]{Liu2025}.
    \begin{theorem}[\cite{Liu2025}]
        Let $d\geq 1$, $n\ge3$ and $\beta>0$. For any $x\in \boxn$,
        \begin{equation}\label{eq: two point functions}
            T_{\beta}^{(n)}(x) - \tau_{\beta}^{\torusn}(x) \leq 2\left(T_{\beta}^{(n)}\star\left(T_{\beta}^{(n)}\right)^{2}\left[T_{\beta}^{(n)} - \tau_{\beta}\right]\star T_{\beta}^{(n)}\right)(x).
        \end{equation}
    \end{theorem}
    We use this theorem to prove~\eqref{eq: chi per bound}.

    \begin{proof}[Proof of \eqref{eq: chi per bound}]
     The lower bound in~\eqref{eq: chi per bound} was established in Proposition~\ref{prop: susceptibility bounds}, so we turn our attention to the remaining upper bound.
        
        Firstly, by the translation invariance of the torus, we have
        \begin{equation}\label{eq: sum conv}
            \sum_{x\in \boxn}(f\star g)(x) = \left(\sum_{x\in \boxn}f(x)\right)\left(\sum_{x\in \boxn}g(x)\right).
        \end{equation}
        Write $T = T_{\beta}^{(n)}$ and $D = T_{\beta}^{(n)} - \tau_{\beta}$. Then applying~\eqref{eq: two point functions} and summing over all $x \in \boxn$, we have by~\eqref{eq: sum conv}
        \begin{equation}\label{eq: chi upper bound}
            \begin{split}
                \chi_{\ZZ^{d},\beta} - \chi_{\boxn, \beta}^{\mathrm{per}} &\leq 2\sum_{x\in \boxn}\left(T\star(T^{2}D)\star T\right)(x),\\
                &=2\left(\sum_{x\in \boxn}T(x)\right)^{2}\sum_{x\in \boxn} T(x)^2D(x),\\
                &\leq2(\chi_{\ZZ^{d},\beta})^{2}\|T\|_{\infty}^{2}\sum_{x\in \boxn}D(x).
            \end{split}
        \end{equation}
        Now applying~\eqref{eq: exp decay tp} with $m = m_{\beta, d}$ gives
        \begin{equation}\label{eq: D L1}
            \begin{split}
                \sum_{x\in \boxn}D(x) &= \sum_{x\in \boxn}\sum_{u \neq 0}\tau_{\beta}(x + nu)\\
                &= \sum_{x\in\ZZ^{d}\setminus \boxn}\tau_{\beta}(x)\\
                &= O(e^{-m n/2}).
            \end{split}
        \end{equation}
        Additionally, by writing $T = \tau_{\beta} + D$, we obtain
        \begin{equation}\label{eq: T inf}
        \begin{split}
            \|T\|_{\infty} &\leq \|\tau\|_{\infty} + \|D\|_{\infty}\\
            &\leq \tau(0) + \sum_{x\in \ZZ^{d}\setminus \Lambda_{n}}\tau_{\beta}(x),\\
            &= 1+ O(e^{-mn/2}).
        \end{split}
        \end{equation}
        Combining~\eqref{eq: chi upper bound}, \eqref{eq: D L1} and~\eqref{eq: T inf} we find there exists a constant $C' = C'(\beta, d)$ such that
        \begin{equation}
            \chi_{\ZZ^{d},\beta} - \chi_{\boxn, \beta}^{\mathrm{per}} \leq C'(\chi_{\ZZ^{d}, \beta})^{2}e^{-mn/2}.
        \end{equation}
        Now the last inequality in~\eqref{eq: chi per bound} easily follows by choosing $k^{\mathrm{per}} = m/2$ and $K^{\mathrm{per}} = C' (\chi_{\ZZ^{d}, \beta})^2$.
    \end{proof}

    We now turn our attention to the case of free boundary conditions and~\eqref{eq: chi 0 bound}. For this we will leverage the Simon--Lieb inequality~\cite{Simon1980, Lieb1980} for the nearest-neighbour Ising model on $\ZZ^{d}$. For a given set $S\subset \ZZ^{d}$, we let $\partial^{-}S:=\{x\in S: \exists y\in \ZZ^{d}\setminus S \text{ with } xy\in E(\ZZ^{d})\}$ be the inner boundary of $S$. Recall from~\eqref{eq: EB def} that for $S \subset \ZZ^{d}$, $E_{\BB}(S)=\{xy: x, y\in S, x-y=\pm e_{i}, i\in [d]\}$ denotes the set of nearest-neighbour edges induced by the vertices in $S$. 
    \begin{theorem}[Simon--Lieb inequality]\label{theorem: simon lieb}
        Let $S$ be a finite subset of $\ZZ^{d}$. Then for all $\beta\leq\beta_{c}$, $x\in S$ and $y \in \ZZ^{d}$
        \begin{equation}\label{eq: simon lieb}
            \avg{\sigma(x)\sigma(y)}_{\ZZ^{d}, \beta} \leq \avg{\sigma(x)\sigma(y)}_{S, \beta}^{0}\ind(y \in S) + \sum_{\substack{u \in S\\ v\in \ZZ^{d}\setminus S \\ uv\in E(\ZZ^{d})}}\avg{\sigma(x)\sigma(u)}_{S, \beta}^{0}\tanh(\beta)\avg{\sigma(v)\sigma(y)}_{\ZZ^{d}, \beta}.
        \end{equation}
    \end{theorem}
    \begin{remark}
        Equation~\eqref{eq: simon lieb} is typically presented with $x \in S$ and $y \in \ZZ^{d}\setminus S$, in which case the first term on the right-hand side of~\eqref{eq: simon lieb} is 0. The slightly more general version above can be proved using the backbone representation for random currents in a manner similar to~\cite[Lemma 2.7]{DC2016}, by decomposing the backbones into those which remain entirely in $S$ and those which do not. 
    \end{remark}

    \begin{proof}[Proof of \eqref{eq: chi 0 bound}]
        Using the definition of the finite-volume and infinite-volume susceptibilities, and then applying translation invariance of the infinite-volume two-point function, we have
        \begin{equation}\label{eq: chi difference 0}
            \begin{split}
                \chi_{\ZZ^{d}, \beta} - \chi_{\boxn, \beta}^{0}
                &=\sum_{y \in \ZZ^{d}}\avg{\sigma(0)\sigma(y)}_{\ZZ^{d}, \beta} - \frac{1}{n^{d}}\sum_{\substack{x, y \in \boxn}}\avg{\sigma(x)\sigma(y)}_{\boxn, \beta}^{0}\\
                &=\frac{1}{n^{d}}\sum_{\substack{x \in \boxn\\y \in \ZZ^{d}}}\left(\avg{\sigma(x)\sigma(y)}_{\ZZ^{d}, \beta} - \avg{\sigma(x)\sigma(y)}_{\boxn, \beta}^{0}\ind(y \in \boxn)\right).
            \end{split}
          \end{equation}
        Note also that Griffiths' inequality~\cite{griffiths1967correlations} implies
        \begin{equation}
          \label{eq:Griffiths}
            \avg{\sigma(x)\sigma(y)}_{\ZZ^{d}, \beta}\geq \avg{\sigma(x)\sigma(y)}_{\boxn, \beta}^{0}, \quad\quad x, y \in \boxn.
        \end{equation}
        Applying the Simon-Lieb inequality~\eqref{eq: simon lieb} with $S = \boxn$, followed by~\eqref{eq:Griffiths} then yields
        \begin{equation}\label{eq: chi diff free bound}
            \begin{split}
                \chi_{\ZZ^{d}, \beta} - \chi_{\boxn, \beta}^{0} &\leq \frac{1}{n^{d}}\sum_{\substack{x \in \boxn\\y \in \ZZ^{d}}}\sum_{\substack{u \in \boxn\\ v\notin \boxn \\ u\sim v}}\avg{\sigma(x)\sigma(u)}_{\boxn, \beta}^{0}\tanh(\beta)\avg{\sigma(v)\sigma(y)}_{\ZZ^{d}, \beta},\\
                &=\frac{\chi_{\ZZ^{d}, \beta}\tanh(\beta)}{n^{d}}\sum_{x\in \boxn}\sum_{\substack{u \in \Lambda_n\\ v\notin \Lambda_n \\ u\sim v}}\avg{\sigma(x)\sigma(u)}_{\Lambda_n, \beta}^{0},\\
                &\leq \frac{2d\,\chi_{\ZZ^{d}, \beta}\tanh(\beta)}{n^{d}}\sum_{x\in \boxn}\sum_{\substack{u \in \partial^{-}\Lambda_n}}\avg{\sigma(x)\sigma(u)}_{\Lambda_n, \beta}^{0},\\
                & \leq \frac{2d\,\chi_{\ZZ^{d}, \beta}^{2}\tanh(\beta)}{n^{d}}|\partial^{-}\boxn|,\\
                & \leq \frac{K^{0}}{n},
            \end{split}
        \end{equation}
        for some $K^{0} = K^{0}(\beta, d)$ sufficiently large. This establishes the upper bound in~\eqref{eq: chi 0 bound}.

        To obtain the lower bound combine~\eqref{eq: chi difference 0}, \eqref{eq:Griffiths} and the lower bound $\avg{\sigma(x)\sigma(y)}_{\ZZ^{d}, \beta}\geq0$, to obtain
        \begin{equation}\label{eq: prelim lower bound 0}
            \begin{split}
                \chi_{\ZZ^{d}, \beta} - \chi_{\boxn, \beta}^{0} &=  \frac{1}{n^{d}}\sum_{\substack{x \in \boxn\\y \in \ZZ^{d}}}\left(\avg{\sigma(x)\sigma(y)}_{\ZZ^{d}, \beta} - \avg{\sigma(x)\sigma(y)}_{\boxn, \beta}^{0}\ind(y \in \boxn)\right),\\
                &\geq \frac{1}{n^{d}}\sum_{\substack{x \in \boxn\\y \in \ZZ^{d}\setminus \boxn}}\avg{\sigma(x)\sigma(y)}_{\ZZ^{d}, \beta},\\
                &\geq \frac{1}{n^{d}}\sum_{\substack{x \in \partial^{-}\boxn\\y \in \ZZ^{d}\setminus \boxn\\x\sim y}}\avg{\sigma(x)\sigma(y)}_{\ZZ^{d}, \beta}
            \end{split}
        \end{equation}
        Now by translation invariance and the $\ZZ^{d}$-symmetry of the infinite-volume two-point function, we have for every $y \sim x$
        \begin{equation}\label{eq: twopoint function symmetry}
            \avg{\sigma(x)\sigma(y)}_{\ZZ^{d}, \beta} = \avg{\sigma(0)\sigma(y-x)}_{\ZZ^{d}, \beta} = \avg{\sigma(0)\sigma(e_1)}_{\ZZ^{d}, \beta}.
        \end{equation}
        Again by Griffiths' inequality~\cite{griffiths1967correlations}
        \begin{equation}\label{eq: two point function lower bound}
            \avg{\sigma(0)\sigma(e_1)}_{\ZZ^{d}, \beta} \geq \avg{\sigma(0)\sigma(e_1)}_{\{0, e_{1}\}, \beta}^{0} = \tanh(\beta),
        \end{equation}
        where $\avg{\, \cdot\,}_{\{0, e_{1}\}, \beta}^{0}$ denotes the expectation with respect to the Gibbs measure defined on the graph with vertex set $\{0, e_{1}\}$ joined by a single edge. Combining~\eqref{eq: prelim lower bound 0}, \eqref{eq: twopoint function symmetry} and~\eqref{eq: two point function lower bound} we obtain
        \begin{equation}
            \begin{split}
                \chi_{\ZZ^{d}, \beta} - \chi_{\boxn, \beta}^{0} &\geq\frac{\tanh(\beta)}{n^{d}}|\partial^{-}\boxn|\\
                &\geq\frac{\tanh(\beta)}{n} := \frac{k^{0}}{n},
            \end{split}
        \end{equation}
        which is precisely the lower bound in~\eqref{eq: chi 0 bound}. 
    \end{proof}

    Finally we establish~\eqref{eq: chi inf bound} using a similar argument.
    \begin{proof}[Proof of~\eqref{eq: chi inf bound}]
        Using translation invariance of the infinite-volume Gibbs measure and then applying the Simon-Lieb and Griffiths inequalities in an identical fashion to Equation~\eqref{eq: chi diff free bound}, we have
        \begin{equation}
            \begin{split}
                0 &\leq \chi_{\ZZ^{d}, \beta} - \chi_{\boxn, \beta}^{\infty},\\
                &=\frac{1}{n^{d}}\sum_{\substack{x\in \boxn\\y\notin \boxn}}\avg{\sigma(x)\sigma(y)}_{\ZZ^{d}, \beta}\\
                &\leq \frac{1}{n^{d}}\sum_{\substack{x \in \boxn\\y \notin \boxn}}\sum_{\substack{u \in \boxn\\ v\notin \boxn \\ u\sim v}}\avg{\sigma(x)\sigma(u)}_{\boxn, \beta}^{0}\tanh(\beta)\avg{\sigma(v)\sigma(y)}_{\ZZ^{d}, \beta},\\
                & \leq \frac{2d\,\chi_{\ZZ^{d}, \beta}^{2}\tanh(\beta)}{n^{d}}|\partial^{-}\boxn|,\\
                & \leq \frac{K^{\infty}}{n},
            \end{split}
        \end{equation}
        for an appropriately chosen $K^{\infty}$. For the lower bound we apply $\avg{\sigma(x)\sigma(y)}_{\ZZ^{d}, \beta}\geq0$ and~\eqref{eq: two point function lower bound} to get
        \begin{equation}
            \begin{split}
                \chi_{\ZZ^{d}, \beta} - \chi_{\boxn, \beta}^{\infty}&\geq\frac{1}{n^{d}}\sum_{\substack{x \in \partial^{-}\boxn\\y \notin \boxn}}\avg{\sigma(x)\sigma(y)}_{\ZZ^{d}, \beta}\\
                &\geq \frac{\tanh(\beta)}{n}=:\frac{k^{\infty}}{n},
            \end{split}
        \end{equation}
        which completes the proof.
    \end{proof}

\backmatter
 \bmhead{Acknowledgements}
 The authors would like to thank Stephen Muirhead for discussions in relation to Section~\ref{sec: proof of prop}. The second author also acknowledges the use of OpenAI's ChatGPT (GPT-5.1 and GPT-5.2) in the early stages of this project, which assisted in developing some aspects of the proof strategy for Theorem~\ref{Theorem: main result 1}, primarily in Section~\ref{sec: proof of theorem}, and motivated Remark~\ref{remark: stein coupling}.
 This research was supported by the Australian Research Council's Discovery Projects funding scheme (Project No. DP230102209). It was also supported through an Australian Government Research Training Program Scholarship [DOI: https://doi.org/10.82133/C42F-K220].
  
\bmhead{Data Availability} Data sharing is not applicable to this article, as no datasets were generated or analysed during the current study.

\bmhead{Conflicts of interest} The authors have no competing interests to declare that are relevant to the content of this article.

\begin{appendices}
    \section{Proof of Proposition~\ref{prop: decay of semigroup}}\label{sec: decay of semigroup}
    In this appendix we give a detailed proof of Proposition~\ref{prop: decay of semigroup}. We recall that for a set $S\subset \ZZ^{d}$, $\partial^{-}S:=\{x\in S:xy\in E(\ZZ^{d}), y\notin S\}$ denotes the inner boundary of $S$. We additionally define $\partial^{+}S:=\{x\notin S: xy\in E(\ZZ^{d}), y\in S\}$ to be the exterior boundary.

    \begin{proof}[Proof of Proposition~\ref{prop: decay of semigroup}]
      In the case of periodic boundary conditions, Equation~\eqref{eq: exp contraction} was established in~\cite[Lemma 2.1]{Lubetzky2016}. It therefore suffices to establish~\eqref{eq: exp contraction} for all $\xi\in \{-1, 0, 1\}^{\ZZ^{d}}$.
      The strategy is to apply~\cite[Theorem 3.1]{Yoshida1997}. Since no assumption of translation invariance of the Hamiltonian is made in~\cite{Yoshida1997}, we may couple the Ising models defined in~\eqref{Eq: gibbs measure} for all $n\geq 3$ and all $\xi\in\{-1, 0, 1\}^{\ZZ^{d}}$ using a single interaction on $\ZZ^{d}$. A log-Sobolev inequality for this enlarged system will then imply the result. Note that for each $n\geq 3$, due to the nearest-neighbour interactions in~\eqref{Eq: hamiltonian}, boundary conditions are specified by a choice of $\xi\in \{-1, 0, 1\}^{\partial^{+}\boxn}$.

      Let $(n_{k}, \xi^{(k)})_{k\geq 0}$ denote an arbitrary enumeration of all pairs $(n,\xi)$ with $n\ge3$ an integer and $\xi\in\{-1,0,1\}^{\partial^{+}_{\boxn}}$. For $k\ge0$, define
        \begin{equation}
            b_{k} := \left(\sum_{\ell=0}^{k-1}n_{\ell}+\lfloor n_{k}/2\rfloor\right)\, e_1, \quad B_{k}:=\Lambda_{n_{k}}+b_{k}.
        \end{equation}
        If $x\in B_k$, then $\sum_{\ell=0}^{k-1}n_{\ell} \le x_1 \le \sum_{\ell=0}^{k}n_{\ell}-1$, which implies the boxes $B_{k}$ are disjoint.

        For each $x, y\in \ZZ^{d}$, define
        \begin{equation}
            J_{xy}=\beta \ind(xy\in E(\ZZ^{d}) \text{ and } x, y\in B_{k} \text{ for some } k\geq 0),
        \end{equation}
        and for $x\in\ZZ^d$ define $h_x$ so that $h_{x} = 0$ on $\ZZ^{d}\setminus\cup_{k\geq0}B_{k}$, and  
        \begin{equation}
            h_{b_{k}+x} = \beta \sum_{\substack{y\in \partial^{+}\Lambda_{n_{k}}\\xy\in E(\ZZ^{d})}}\xi^{(k)}_{y}, \qquad x\in \Lambda_{n_{k}}, \, k\ge0.
        \end{equation}
        Now define the interaction $\Phi=\{\Phi_X\}_{X\Subset\ZZ^d}$ so that for all distinct $x,y\in\ZZ^d$
        \begin{equation}
            \Phi_{\{x\}}(\sigma):= -h_{x}\sigma(x), \quad \Phi_{\{x, y\}}(\sigma) = -J_{x, y}\sigma(x) \sigma(y), \qquad \sigma\in\{-1,1\}^{\ZZ^d},
        \end{equation}
        with all other $\Phi_X=0$. Since these interactions are uniformly bounded and have range at most 1, assumptions $(\Phi\text{-}1)$, $(\Phi\text{-}2)$ in~\cite{Yoshida1997} are satisfied. 

        For any finite nonempty $V\subset \ZZ^d$ and field $\mathfrak{h}\in\RR^V$ define
        \begin{equation}\label{eq: gibbs V}
          \mu_{V,J,\mathfrak{h}}(\sigma) :=\frac{1}{Z_{V,J,\mathfrak{h}}}
          \exp\left(\sum_{xy\in E_{\BB}(V)} \,J_{xy}\,\sigma(x)\sigma(y)
            + \sum_{x\in V}\mathfrak{h}_{x}\sigma(x)\right), \quad \sigma\in \{-1, 1\}^{V}.
        \end{equation}
        The Gibbs measure on $V$ with interaction $\Phi$ and boundary condition $\omega$ then corresponds to $\mu_{V,J,\eta}$ with
        the field $\eta\in\RR^V$ defined by
        \begin{equation}
          \eta_x:=\left(h_x+\sum_{y\notin V}J_{xy}\omega(y)\right).
        \end{equation}

        We also consider the corresponding heat-bath process with stationary measure $\mu_{V,J,\eta}$, which has flip rates $c^{V,J,\eta}_{x}(\sigma):=\frac{1}{2}[1-\sigma(x)\tanh(S_{x}^{V,J,\eta}(\sigma))]$ where   
        \begin{equation}
            S_{x}^{V,J,\eta}(\sigma):= \sum_{\substack{y\in V\\xy\in E_{\BB}(V)}}J_{xy}\sigma(y) +\eta_{x}.
          \end{equation}
       We denote the corresponding semigroup by $P_{t}^{V,J,\eta}$.
       Since $|S_{x}^{V,J,\eta}(\sigma)|\leq 2d\beta$, the rates satisfy
        \begin{equation}
            0<\frac{1}{2}[1-\tanh(2d\beta)]\leq c_{x}^{V,J,\eta}(\sigma)\leq \frac{1}{2}[1+\tanh(2d\beta)].
        \end{equation}
        Together with the facts that $c^{V,J,\eta}_{x}$ is a local function and the dynamics are reversible, we conclude that hypotheses $(\text{R-}1)$, $(\text{R-}2)$ and $(\text{R-}3)$ from~\cite{Yoshida1997} are satisfied. It follows from~\cite[Theorem 3.1]{Yoshida1997} that provided $P_t^{V,J,\eta}$ satisfies a log-Sobolev inequality, uniformly in $V$ and $\omega$, we will obtain a uniform exponential relaxation estimate for any local function, from which~\eqref{eq: exp contraction} will follow as a corollary.
 
        Fortunately, such a log-Sobolev inequality was recently given in~\cite{BauerschmidtDagallier2024}.
        Note that $\mu_{V,J,0}$ corresponds to the zero-field Ising measure on $V$ with free boundary conditions and couplings $(J_{xy})_{x, y\in \ZZ^{d}}$. Let $\<\,\cdot\,\>_{V,J,0}$ denote expectation with respect to $\mu_{V,J,0}$.
        It follows by Griffiths' inequality~\cite{griffiths1967correlations} that for all $\beta<\beta_{c}(d)$
        \begin{equation}
          \sup_{x\in V}\sum_{y\in V}\<\sigma(x)\sigma(y)\>_{V,J,0}
          \leq \sum_{y\in \ZZ^{d}}\<\sigma(0)\sigma(y)\>_{\ZZ^{d}, \beta} = \chi_{\ZZ^{d}, \beta}<\infty.
        \end{equation}
        We may then appeal to~\cite[Theorem 1.1]{BauerschmidtDagallier2024}, which provides a bound for the log-Sobolev constant of $\mu_{V,J,\eta}$ which is uniform in $V$ and $\eta$, and therefore also in $\xi$ and $\omega$.
        Thus \cite[Equation (3.1)]{Yoshida1997} is satisfied, and so applying~\cite[Theorem 3.1 (1c)]{Yoshida1997} to the function $f_{x}(\sigma) = \sigma(x)$ yields 
        \begin{equation}
          \label{eq: Phi exponential relaxation}
       \sup_{\substack{V, \omega, \sigma}}|P_{t}^{V, J,\eta}f_{x} - \<\sigma_{x}\>_{V,J,\eta}|\leq C_{\beta, d}e^{-c_{\beta, d}t}.
     \end{equation}
The constants $c_{\beta,d}$ and $C_{\beta,d}$ are independent of $x$ by the admissibility condition~\cite[Equation (2.17)]{Yoshida1997}, and are uniform in $n$ and $\xi$, since the interaction $\Phi$ is fixed and the estimate is uniform in $V$ and $\omega$.

Finally, we note that, by construction, when $V=B_k$ we have $\eta_x=h_x$, for all $x\in B_k$. 
        Moreover, for any integer $n\ge3$ and $\xi\in\{-1,0,1\}^{\ZZ^d}$ there exists $k\ge0$ such that $n_k=n$,
        $\xi^{(k)}=\xi|_{\partial^+\boxn}$, and for all $\sigma\in\{-1,1\}^{\boxn}$ 
        \begin{equation}
          \label{eq: translated gibbs measure}
          \rho_{\boxn,\beta}^{\xi}(\sigma) = \mu_{B_k, J, h}(\tau_k\sigma)
        \end{equation}
        where $\tau_k$ is the translation defined by $(\tau_k\sigma)(x):=\sigma(x-b_k)$.
        Similarly, it is easily verified that $c_{x+b_k}^{B_k,J,h}(\tau_k\sigma)=c_x^{\xi}(\sigma)$ for all $x\in\boxn$ and $\sigma\in\Omega^{\boxn}$, from which it follows that
        \begin{equation}
          \label{eq: translated semigroup}
          \EE [X_t^{\xi,\sigma}(x)] = P_t^{B_k,J,h} f_{x+b_k}(\tau_k\sigma) 
        \end{equation}
        where $f_z(\nu):=\nu(z)$ for all $z\in B_k$ and $\nu\in\Omega^{B_k}$.

        Restricting~\eqref{eq: Phi exponential relaxation} to $V\in \{B_{k}\}_{k\geq0}$, and using~\eqref{eq: translated semigroup} and~\eqref{eq: translated gibbs measure} then implies that
        \begin{equation}
            \sup_{\sigma, x, \xi}\left|\EE[X_t^{\xi,\sigma}(x)] - \avgl{\sigma(x)}\right| \leq C_{\beta,d}\,e^{-c_{\beta, d}t}
        \end{equation}
        where $C_{\beta, d}$ and $c_{\beta, d}$ are independent of $n$ and $\xi$.
    \end{proof}
\end{appendices}


\end{document}